\documentclass[11pt,a4paper]{article}
\usepackage{amsmath}
\usepackage{amssymb}
\usepackage{amsthm}
\usepackage{graphicx}
\usepackage{xcolor}
\usepackage{bm}
\usepackage{etoolbox}
\usepackage{xspace}
\usepackage{here}
\usepackage{color}
\usepackage[runin]{abstract}
\usepackage{appendix}
\numberwithin{equation}{section}

\renewcommand\v[1]{{\bm {#1}}}
\newcommand\inner[2]{{\left\langle{#1},{#2}\right\rangle}}
\newcommand\Z{{\mathbb Z}}
\newcommand\R{{\mathbb R}}

\newcommand{\AltQuad}{alternating quadruple\xspace}
\newcommand{\AltQuads}{alternating quadruples\xspace}
\newcommand{\vfunc}{\vec}
\theoremstyle{definition}
\newtheorem{proposition}{Proposition}[section]
\newtheorem{example}[proposition]{Example}
\newtheorem{theorem}[proposition]{Theorem}
\newtheorem{remark}[proposition]{Remark}
\newtheorem{convention}[proposition]{Convention}

\newtheorem{lemma}[proposition]{Lemma}
\newtheorem{definition}[proposition]{Definition}

\newtheorem{mainresult}[proposition]{Main Result}
\newtheorem*{theorem*}{Theorem}
\numberwithin{figure}{section}
\numberwithin{table}{section}
\newcounter{claimcounter}
\AtBeginEnvironment{proof}{\setcounter{claimcounter}{0}}
\newenvironment{claim}[1][]{\refstepcounter{claimcounter}\par\medskip\noindent\textbf{Claim \theclaimcounter. #1}\itshape}{\par}
\newenvironment{claimproof}[1][]{\medskip\noindent\textit{Proof of Claim~\theclaimcounter. #1}\ignorespaces}{%
  \hfill$\blacksquare$%
  \par\medskip
}
\begin{document}
\title{
  \LARGE
  \bfseries 
  Stable configurations of entangled systems with repulsive interactions
}
\date{}
\author{
  Motoko Kotani\footnote{The Advanced Institute for Materials Research (AIMR), Tohoku University, JAPAN}, \mbox{} 
  Hisashi Naito\footnote{Graduate School of Mathematics, Nagoya University, JAPAN}, \mbox{}
  Naoki Sakata\footnote{The Advanced Institute for Materials Research (AIMR), Tohoku University, JAPAN}, \mbox{}
  and 
  Eriko Shinkawa\footnote{Mathematical Science Center for Co-creative Society, Tohoku University, JAPAN, and 
    The Advanced Institute for Materials Research (AIMR), Tohoku University, JAPAN}
}
\maketitle
\begin{abstract}
  Entangled systems are prevalent in both biological and synthetic materials.
  This study examines the stable configurations of weaves consisting of two families of intertwined threads, such as warp and weft threads.
  By analyzing the steepest descent flow of an energy functional featuring repulsive interactions, we develop a framework for identifying stable states in $\R^3$.
  Although a weave consists of one-dimensional threads that do not intersect each other, it behaves collectively like a two-dimensional object.
  To describe this phenomenon, we define a non-separable component of a weave as a ``layer'' and establish the existence and uniqueness of its stable configuration.
  Furthermore, we show that two distinct layers drift apart with an asymptotic growth rate of $t^{1/3}$ as $t \to \infty$.
\end{abstract}

\section{Introduction}
\label{sec:introduction}
Entanglement is a topological object widely observed in filamentary and networked systems, ranging from biological tissues to synthetic materials.
In nature, interwoven fiber networks play a key role in the mechanical stability of biological tissues, such as plant cell walls and connective tissues. They have also been used for a long time in the form of fibers and textiles in industry.
The field of synthetic chemistry has seen substantial progress in recent years, particularly in designing materials with molecularly controlled interactions and network structures.
For example, mechanically interlocked architectures, such as molecular weaves, as well as interpenetrating frameworks in metal-organic frameworks (MOFs) and covalent organic frameworks (COFs), are receiving increasing attention.
Furthermore, researchers are developing fabrics made from carbon nanotube yarns and fibers with high-performance physical properties.

In theoretical studies of entanglement and network structures, geometric descriptions and systematic frameworks have been developed.
Liu--O'Keeffe--Treacy--Yaghi proposed a systematic geometrical study of periodic networks and weavings as a library for reticular chemistry~\cite{lkty}.
Hyde--Chen--O'Keeffe applied topological weaving descriptions to relate planar layers with 3D frameworks through geometric warping~\cite{okeeffe}.
In this context, Evans--Hyde and co-workers have intensively studied periodic entanglements of closed components and nets, providing systematic geometrical descriptions of entangled structures (see~\cite{evans-hyde0, evans-hyde1, evans-hyde2, evans-hyde3} and references therein).

Textile structures are studied using the tools of knot theory (see, e.g., \cite{grishanov1, grishanov2, grishanov3, kawauchi}).
For instance, textile structures can be mathematically described as patterns of doubly periodic link diagrams (see, e.g.,~\cite{diamantis1}).
To classify these structures, knot invariants applicable to periodic diagrams have been investigated.
In particular, partial classification results have been obtained, and it has been proven that certain fundamental patterns are distinct (see, e.g.,~\cite{diamantis2, grishanov2}).

In this paper, an ``entangled'' system refers to a disjoint union of connected spatial graphs.
While the term ``entangled'' generally implies that the components are non-separable in space, our framework also allows for configurations where they can be separated.
Intuitively, we refer to a maximal set of such mutually non-separable graphs as a ``layer''. 
Consequently, an entangled system is formulated as a disjoint union of these layers. 
A precise mathematical definition of a layer is provided later in Definition~\ref{definition:weaving:weaving_layer}.
From the viewpoint of topology, we regard these entangled structures as equivalent under ambient isotopy in the three-dimensional space $\mathbb{R}^3$.
Namely, they can be continuously deformed as long as their topological type and the underlying combinatorial structure as a graph are preserved.

In contrast, even if entangled structures represent the same topological type, they exhibit different functions when their physical configurations in $\R^3$ are different.
In many important systems, molecules and materials adopt specific, stable configurations that are dictated by physical laws.
Such stable configurations play a central role in determining the functions and properties of a material.
Therefore, it is natural to ask how one can mathematically describe and identify physically stable configurations for a given entangled structure.

Kotani--Sunada~\cite{Kotani-Sunada} (cf.~\cite{Sunada}) defined a ``crystal lattice'', a topological description of crystal structures in terms of connected, infinite periodic graphs, where atoms are modeled as vertices and bonds as edges.
Furthermore, for a crystal lattice, they proved the existence and uniqueness (up to Euclidean isometries) of its stable configuration and called it the ``standard realization'' of a given crystal lattice.
It is defined as a periodic realization that minimizes an energy functional under the constraint that the volume of the fundamental domain is fixed.
Here, the energy is formulated as an elastic (spring) energy, treating the edges of the periodic graph as springs.

The Kotani--Sunada theory does not apply to entangled systems because the energy depends only on the edges within a single crystal lattice, leaving distinct components without interaction.
To analytically capture entanglement, we extend the Standard Realization with Repulsive Interactions (SRRI) model, which was originally applied to single crystal lattices~\cite{dechant}, to multi-component systems.
In this approach, we formalize each component as an interacting point set.
The internal elastic energy defines its periodic graph structure and controls its internal configuration, while an added repulsive term between vertices of different components controls their relative positions.
Although this formulation does not strictly prohibit geometric intersections, the energetic barriers penalize overlap.
This penalty can reproduce the topological constraints of the system through mechanical conditions.
Our overarching goal is to identify the stable configurations of these multi-component periodic graphs.
In this paper, we focus on textile structures, which serve as simple yet representative model systems in materials science, to achieve this goal.

\begin{figure}[htbp]
  \centering
    \includegraphics[width=0.4\linewidth]{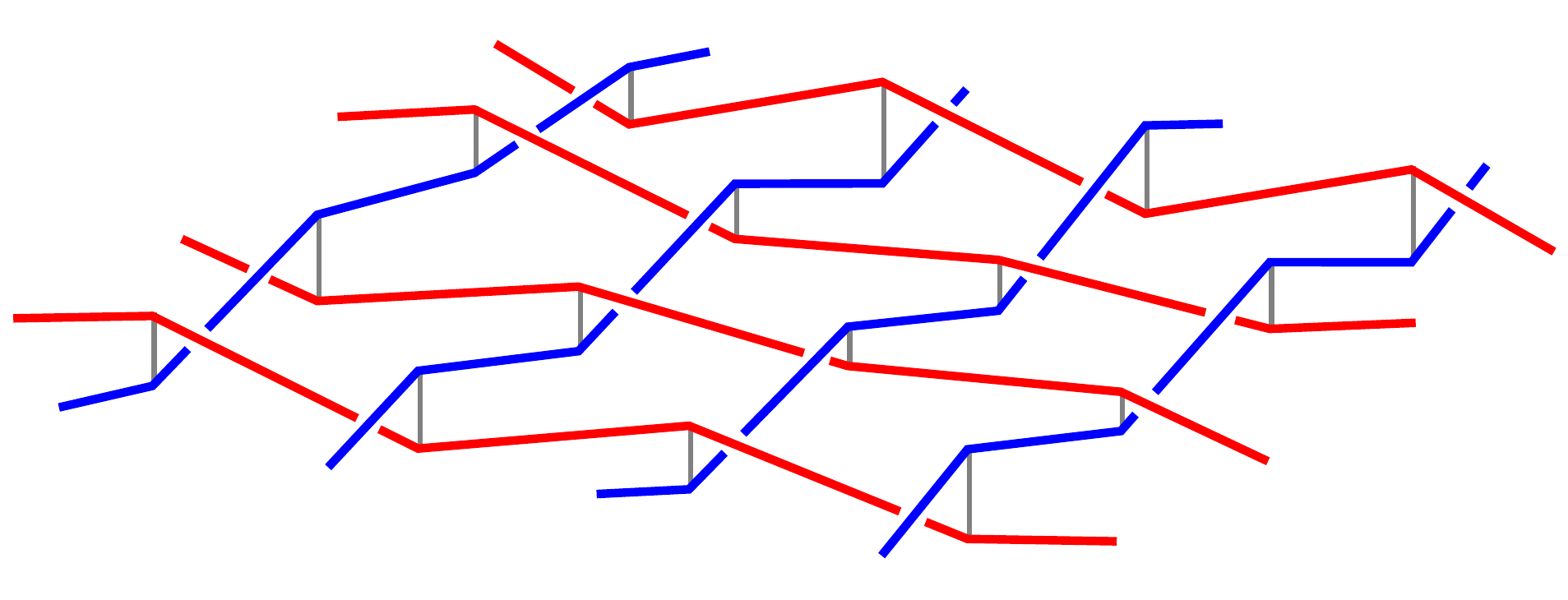}
  \caption{A model of weave studied in this paper. Repulsive forces act between vertices connected by black line segments.}
  \label{fig:intro:models}
\end{figure}

Our model is a mathematical model for textile structures called ``weaves'' (see Figure~\ref{fig:intro:models}), which consist of connected components of one-dimensional manifolds (polygonal paths) referred to as ``threads''.
Although a weave consists of one-dimensional threads that do not intersect each other, it behaves collectively like a two-dimensional object.
To describe this phenomenon, we introduce a notion of ``weave-connectedness'' in Section~\ref{sec:weaving}.

Before stating the main results, we give some terminology.
Roughly speaking, we call the system \emph{separable} if there exists an embedded periodic surface in $\R^3$ which is disjoint from the system and homeomorphic to $\R^2$, dividing the space into two regions such that each region contains at least one component.
The precise definition is given in Definition~\ref{definition:3:separable}.
Otherwise, the system is called \emph{non-separable}.
Intuitively, a ``layer'' of a weave is a maximal non-separable ``subweave'', and is proved to be a weave-connected component, except in special cases (Definition~\ref{definition:weaving:weaving_layer}).

In the SRRI model, a \emph{stable configuration} means a stable equilibrium of the steepest descent flow associated with the SRRI energy
(Definition~\ref{definition:ode:stable}).
The existence of a stable configuration for a weave is characterized by the notion of weave-connectedness.
In particular, we show that the system admits a stable configuration if it consists of a single layer.

The main results of this paper can be summarized as follows.

\begin{mainresult}[Theorem~\ref{claim:weaving:ode:1-component}]
  If the system is non-separable (or consists of a single layer), then the SRRI model admits a unique stable configuration.
\end{mainresult}
\begin{mainresult}[Theorem~\ref{claim:weaving:ode:untangle}]
  If the system is separable, 
  then two ``layers'' with repulsive interactions drift apart on the order of 
  $t^{1/3}$ as time $t \to \infty$.
  Furthermore, each layer converges to a stable configuration of the single-layer system obtained by removing the other layers.
\end{mainresult}

Similar results are obtained for an entangled pair of identical graphs.
The settings and proofs are illustrated with the alternatingly entangled square lattices in Appendix B.

\section{Weaves and weave layers}
\label{sec:weaving}
In this section, we discuss the concept of weaves and their separability within three-dimensional space.
As introduced in Section~\ref{sec:introduction}, separability is a topological notion where a maximal collection of disjoint threads that cannot be spatially separated from one another intuitively forms a single ``layer'' (i.e., a single piece of cloth).
This topological property is translated into the combinatorial structure of a weave through a property we introduce in this paper called~``weave-connectedness''.

\subsection{Weaves and weave-connectedness}
\label{sec:weaving:weave-connectedness}

Various studies on topological weaves have been developed.
For example, joint works involving the first author can be found in \cite{fks1,fks2}.
In this paper, we focus on a special type of weave called a ``multiaxial weave'' 
(which corresponds to an ``untwisted weave'' in \cite{fks1} with translational symmetry).
\par
\begin{figure}[htbp]
  \centering
  \begin{minipage}{0.48\linewidth}
    \centering
    \includegraphics[width=0.8\linewidth,page=1]{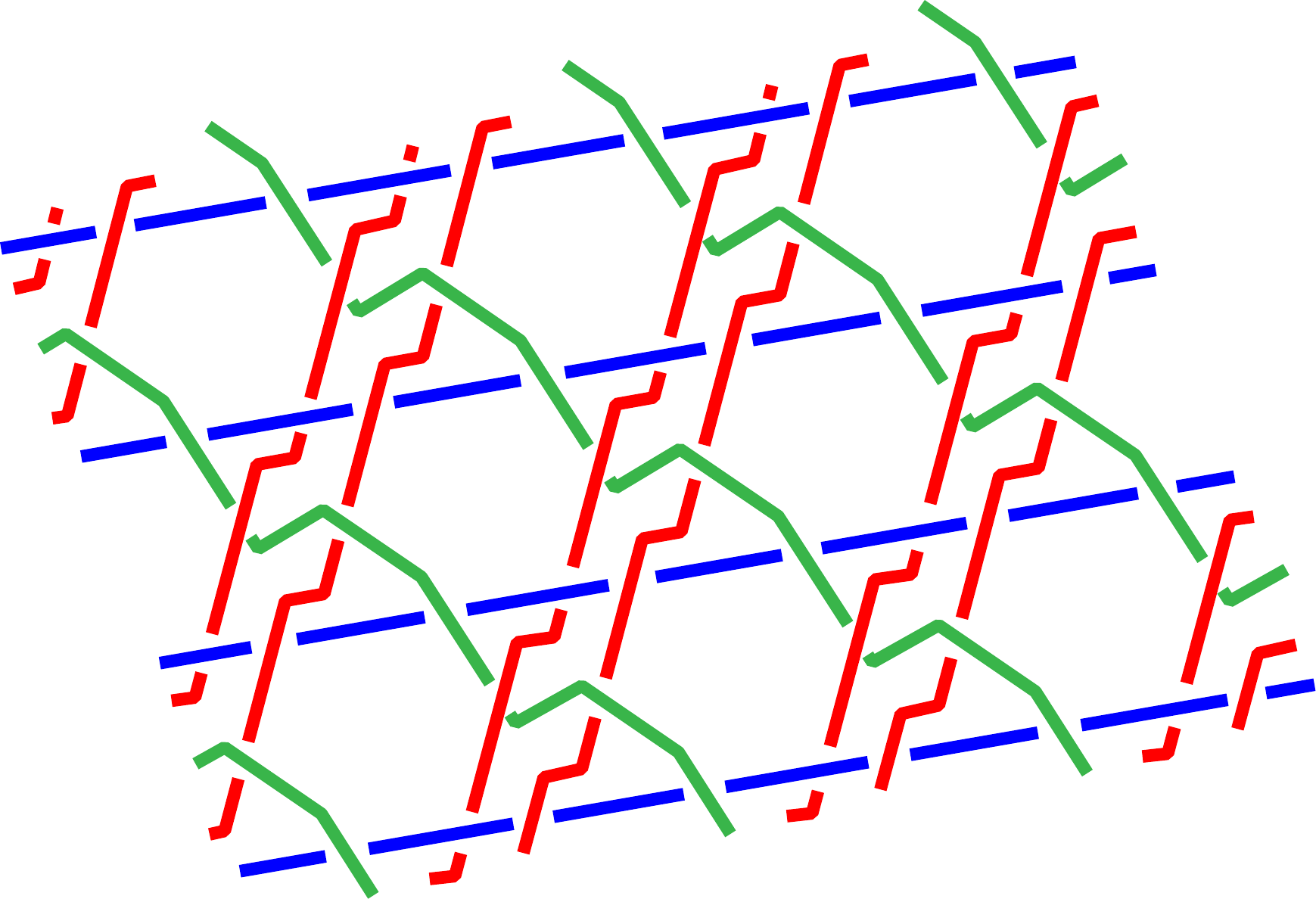}\\
    (a)
  \end{minipage}
  \begin{minipage}{0.48\linewidth}
    \centering
    \includegraphics[width=0.6\linewidth,page=2]{Figure/multiaxial_weave.pdf}\\
    (b)
  \end{minipage}
  \caption{(a) A multiaxial weave. (b) The regular projection of the weave. Each crossing has the over/under information of the threads similar to those in knot diagrams. The weave is preserved under the translations $\tau_1$ and $\tau_2$. The greyed region in (b) shows a motif.}
  \label{fig:2:multiaxial-weave}
\end{figure}

A \emph{multiaxial weave} is defined as a collection of several families that consist of mutually disjoint ``threads''. 
Each \emph{thread} is a connected $1$-dimensional manifold embedded in three-dimensional space $\mathbb{R}^3$.
Furthermore, the orthogonal projection $\pi$ from $\R^3$ to the $xy$-plane maps each family of threads bijectively onto a discrete family of parallel straight lines in the $xy$-plane.
In this paper, we assume that the multiaxial weave is periodic under translations, $\tau_1$ and $\tau_2$, defined by linearly independent vectors in the $xy$-plane (see Figure~\ref{fig:2:multiaxial-weave}).
We call such a translation a \emph{periodic translation}.

By assigning over/under information to each intersection point such that $\pi$ provides a \emph{regular projection} (analogous to knot diagrams in knot theory), the weave can be completely described by a planar diagram. 
Since multiaxial weaves exhibit the assumed periodicity, the entire topological structure can be reconstructed from a fundamental building block called a \emph{motif} (see Figure~\ref{fig:2:multiaxial-weave}(b)). 
This motif can be regarded as a special case of \emph{periodic tangles} (see~\cite{Kotorii-Yoshida} for details).

We focus on a specific class of multiaxial weaves, where the threads are defined as polygonal paths to easily apply the SRRI model.
While the general definition admits arbitrary numbers of thread families, in this paper, we restrict our attention to the case with exactly two families.

\begin{figure}[htbp]
  \centering
  \begin{minipage}[b]{0.49\linewidth}
    \centering
    \includegraphics[width=0.9\linewidth]{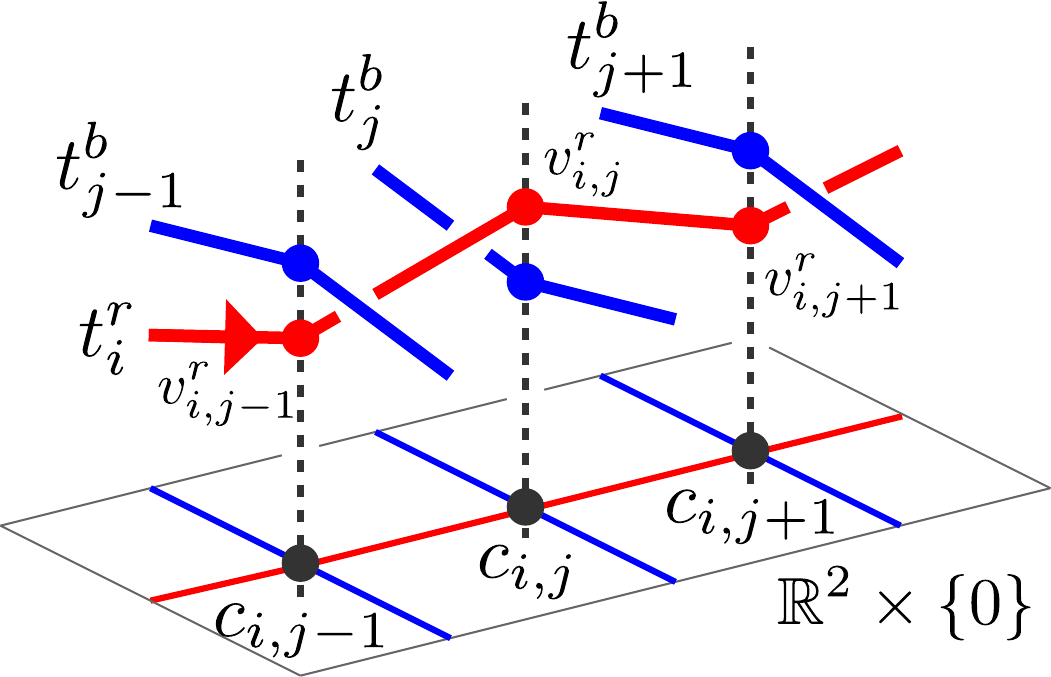}\\
    (a)
  \end{minipage}
  \begin{minipage}[b]{0.30\linewidth}
    \centering
    \includegraphics[width=0.9\linewidth]{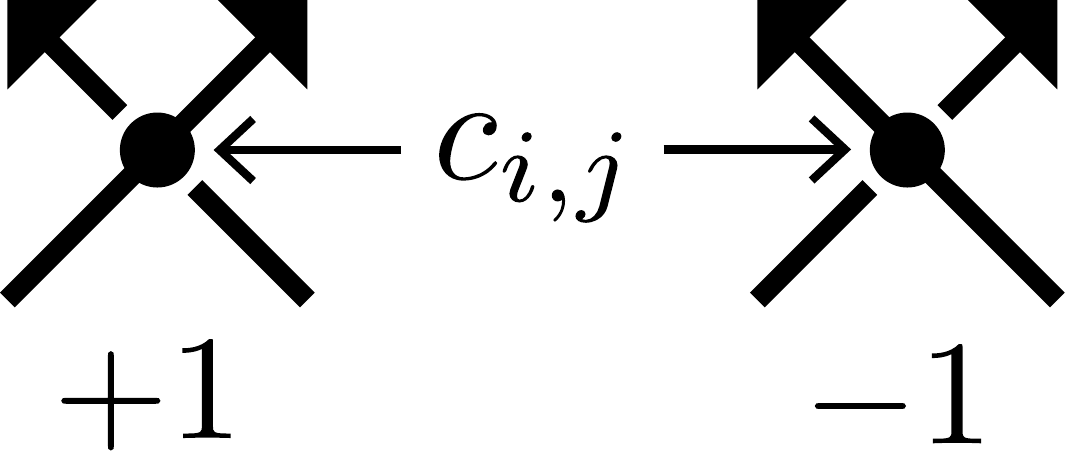}\\
    \vspace{3em}
    (b)
  \end{minipage}
  \caption{(a) A part of threads in a polygonal weave. Dotted lines represent the preimage of the crossings. The order of vertices gives a natural orientation for each thread. (b) The sign $S(i, j)$ at the crossing $c_{i, j}$ of $t^r_i$ and $t^b_j$ with respect to the natural orientations.}
  \label{fig:2:polygonal-weave}
\end{figure}

\begin{definition}
  \label{definition:3:weave}
  A multiaxial weave $W = (T^r, T^b)$ is called a \emph{polygonal weave with two families} (\emph{weave} for short) if the families of threads $T^r = \{t^r_i\}_{i \in \mathbb{Z}}$ and $T^b = \{t^b_j\}_{i \in \mathbb{Z}}$ satisfy the following conditions (see Figure~\ref{fig:2:polygonal-weave}(a)):
  \begin{enumerate}
    \item Let $\pi$ be the orthogonal projection defined by $\pi(x, y, z) = (x, y, 0)$.
      A thread $t^r_i$ (resp. $t^b_j$) has a vertex $v^r_{i, j}$ (resp. $v^b_{i, j}$) in $\pi^{-1}(c_{i, j})$ for each $(i, j) \in \mathbb{Z} \times \mathbb{Z}$.
      Here, $c_{i, j}$ is the crossing of $t^r_i$ and $t^b_j$ in the $xy$-plane, i.e., $c_{i, j} = \pi(t^r_i) \cap \pi(t^b_j)$.
      Furthermore, each of threads $t^r_i$ and $t^b_j$ does not have any vertices except for the preimage of the crossing points.
    \item Each thread $t^r_i$ (resp. $t^b_j$) consists of line segments joining the vertices $v^r_{i, j}$ and $v^r_{i, j+1}$ (resp. $v^b_{i, j}$ and $v^b_{i+1, j}$) for each $j \in \mathbb{Z}$ (resp. $i \in \mathbb{Z}$).
  \end{enumerate}
  For convenience, the threads belonging to the families $T^r$ and $T^b$ are referred to as \emph{red threads} and \emph{blue threads}, respectively, and are depicted in these colors in Figure~\ref{fig:2:polygonal-weave} and subsequent figures.

  A pair $W_1 = (T_1^r, T_1^b)$ is called a \emph{subweave} of $W = (T^r, T^b)$ if there exist index sets $I^r, J^b \subset \Z$ such that the following hold:
  \begin{enumerate}
     \item $T_1^r = \{ t^r_i \in T^r \mid i \in I^r \}$ and $T_1^b = \{ t^b_j \in T^b \mid j \in J^b \}$.
     \item  There exist two periodic translations along linearly independent vectors which preserve $T^r$ and $T^b$ setwise, while also preserving $T^r_1$ and $T^b_1$ setwise.
  \end{enumerate}
Throughout this paper, for a weave $W$, the notation $(\{t^r_{i_1}, \dots, t^r_{i_n} \}, \{ t^b_{j_1}, \dots, t^b_{j_m} \})_W$ implies that the subweave includes not only $t^r_{i_1}, \dots, t^r_{i_n}$, $t^b_{j_1}, \dots, t^b_{j_m}$ but also all threads generated by the translation symmetry of $W$.

  By abuse of notation, a thread $t$ is said to be \emph{in} a weave or subweave $W = (T^r, T^b)$ if $t \in T^r \cup T^b$.
  Furthermore, $t \in W$ is also used to mean $t \in T^r$ or $t \in T^b$.
\end{definition}

Any multiaxial weave can be continuously deformed so that the orthogonal projection $\pi$ becomes a regular projection.
By the definition of a polygonal weave, each thread has a natural orientation defined by the index.
We can obtain the spatial configuration of a weave from its regular projection and the ``sign'' defined as follows.

\begin{definition}
  \label{definition:weave:sign}
  By the definition of polygonal weave, the vertices of a thread $t^r_i$ are ordered as $\dots,\, v^r_{i, j-1},\, v^r_{i, j},\, v^r_{i, j+1},\, \dots$ for each $j \in \Z$.
  Thus, a natural orientation is defined on $t^r_i$ by its indices.
  Similarly, a natural orientation is also defined on $t^b_j$.

  The \emph{sign} $S(i, j) \in \{ +1, -1 \}$ at the crossing $c_{i, j}$ of $t^r_i$ and $t^b_j$ is defined similarly to the crossing sign of an oriented link diagram in knot theory (see Figure~\ref{fig:2:polygonal-weave}(b)).
  We refer to $S(i, j)$ as \emph{the sign between $t^r_i$ and $t^b_j$}.
\end{definition}

Any polygonal weave can be transformed into a more intuitive form through a slight deformation.
Specifically, weaves can be continuously deformed so that the threads within each family are equally spaced.
In general, a periodic translation is not necessarily parallel to the image of a thread under the orthogonal projection $\pi$.
However, a simple argument shows that we can choose a pair of periodic translations to be parallel to the projected lines of two families, respectively.
(Note that the size of the fundamental region may change.)
Thus, we can always deform a weave, applying a reflection if needed, to align its threads with the \emph{square lattice} in the $xy$-plane.

\begin{figure}[htbp]
  \centering
  \begin{minipage}[b]{0.38\linewidth}
    \centering
    \includegraphics[width=1.0\textwidth]{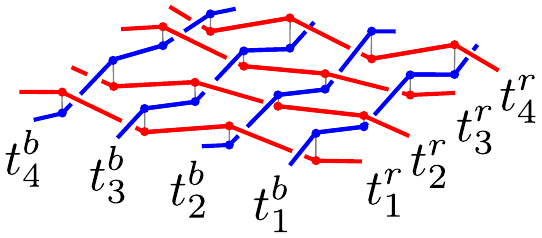}\\
    \vspace{1em}
    (a)
  \end{minipage}
  \begin{minipage}[b]{0.38\linewidth}
    \centering
    \includegraphics[width=0.6\textwidth]{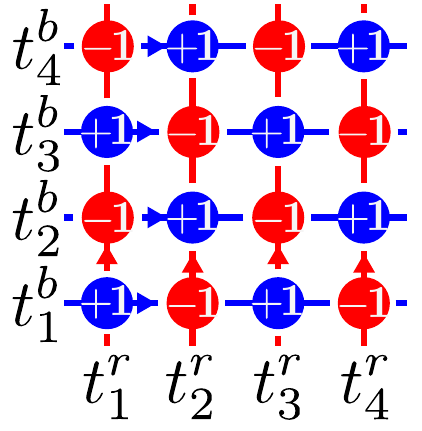}\\
    (b)
  \end{minipage}
  \caption{(a) A weave on the square lattice satisfying Convention~\ref{convention:weave}. (b)~The regular projection of the weave in (a). The values $+1$ or $-1$ shown at each crossing represent its sign.}
  \label{fig:2:square-lattice-weave}
\end{figure}

\begin{convention}
  \label{convention:weave}
  To simplify the proofs, it is convenient to deal with weaves in a standard form.
  Thus, we may assume the following conditions as needed:
  \begin{enumerate}
    \item The vertices $v^r_{i, j}$ and $v^b_{i, j}$ of the threads $t^r_i$ and $t^b_j$ are located at $(i, j, z^r_{i, j})$ and $(i, j, z^b_{i, j})$, respectively.
      So, the projected lines of the threads align with the square lattice in the $xy$-plane, where the images of $T^r$ and $T^b$ consist of lines parallel to the $y$-axis and the  $x$-axis, respectively.
    \item The periodic translations $\tau_1$ and $\tau_2$ are given by vectors of equal length along the $x$- and $y$-axes, respectively.
      Translations along the $x$-axis and $y$-axis with periods $L$ and $M$, respectively, can be made to have the same length by taking $LM$ as their common period.
  \end{enumerate}
  An example of a weave satisfying these conditions is shown in Figure~\ref{fig:2:square-lattice-weave}.
  While the orthogonal projected images of the threads of a subweave are also parallel to the $x$- and $y$-axes, the subweave need not fill the entire square lattice.
\end{convention}

\begin{remark}
  \label{remark:weave:sign}
  The signs of the crossings for a weave under Convention~\ref{convention:weave} can be determined by the relative heights of the threads at the crossings as follows:
  \[
    S(i, j) = 
    \begin{cases} 
      +1 & \text{if } z^b_{i,j} > z^r_{i,j}, \\
      -1 & \text{if } z^b_{i,j} < z^r_{i,j}.
    \end{cases}
  \]
\end{remark}

Next, we will define the notion of ``weave-connectedness''.
When two threads are weave-connected, it is shown that the two threads are contained within the same ``layer'' and cannot be separated by any (topological) plane in the space.

\begin{definition}
  \label{definition:section2:altquad}
  A subweave $(\{t^r_{i_1}, t^r_{i_2}\}, \{t^b_{j_1}, t^b_{j_2}\})_W$ of $W = (T^r, T^b)$ 
  is called an \emph{\AltQuad (of threads)} of $W$, if the sign $S$ satisfies 
  \begin{displaymath}
    (S(i_1, j_1), S(i_2, j_1), S(i_2, j_2), S(i_1, j_2))
    =
    (+1, -1, +1, -1) \text{ or } (-1, +1, -1, +1).
  \end{displaymath}

  A pair of threads $t$ and $t'$ in $W$ are called \emph{weave-connected} if there exists a sequence of \AltQuads $U_0, \dots, U_l$ such that
  \begin{enumerate}
    \item $t \in U_0$ and $t' \in U_l$.
    \item For each $k = 0, \dots, l-1$, $U_k$ and $U_{k+1}$ share a thread.
  \end{enumerate}
  It is clear that weave-connectedness is an equivalence relation on the set of threads in $W$.
  We call an equivalence class of weave-connected threads a \emph{weave-connected component} (see Figure~\ref{fig:weave-connected}).

  A thread is called an {\em untangled thread} if it is not weave-connected to any other thread.
  \begin{figure}[htbp]
    \centering
    \begin{minipage}{0.23\textwidth}
      \centering
        \includegraphics[width=0.95\linewidth]{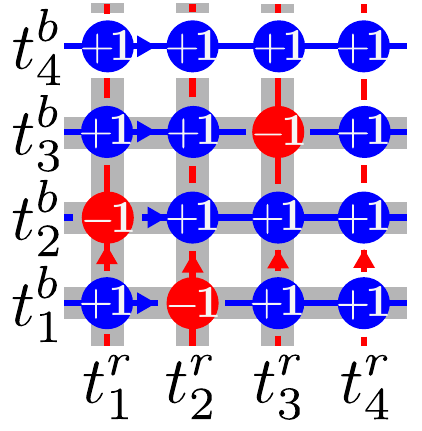}
    \end{minipage}
    \begin{minipage}{0.35\textwidth}
      \begin{align*}
        U_1 &= (\{ t^r_1, t^r_2 \}, \{ t^b_1, t^b_2 \})_W\\
        U_2 &= (\{ t^r_1, t^r_3 \}, \{ t^b_2, t^b_3 \})_W\\
        U_3 &= (\{ t^r_2, t^r_3 \}, \{ t^b_1, t^b_3 \})_W
      \end{align*}
      \vspace{0.4em}
    \end{minipage}
    \caption{An example of weave-connected threads. This weave contains only three \AltQuads $U_1$, $U_2$, and $U_3$. The threads surrounded by the grey area are in the same weave-connected component, $(\{ t^r_1, t^r_2, t^r_3 \}, \{ t^b_1, t^b_2, t^b_3 \})_W$.}
    \label{fig:weave-connected}
  \end{figure}
\end{definition}

\begin{remark}
  By definition, excluding untangled threads, a weave-connected component contains at least one \AltQuad.
\end{remark}

We will give a precise definition of separability for weaves in the following.

\begin{definition}
  \label{definition:3:separable}
  A weave $W$ is said to be \emph{separable} if there exists an embedded surface $\widetilde{\Sigma}$ in $\mathbb{R}^3$ that is disjoint from $W$ and satisfies the following conditions:
   \begin{enumerate}
     \item There exists a pair of periodic translations $\tau_1$ and $\tau_2$ of $W$ such that the translations are defined by linearly independent vectors in the $xy$-plane, and the surface $\widetilde{\Sigma}$ is invariant under these translations.
     \item The quotient, $\Sigma$, of $\widetilde{\Sigma}$ by $\tau_1$ and $\tau_2$ is homeomorphic to $S^1 \times S^1$
   \item The surface $\Sigma$ is parallel to the boundary of the quotient of $\R^2 \times [-B, B]$, where $B$ is a positive number such that the threads of $W$ are contained in $\R^2 \times [-B, B]$.
   \item The surface $\widetilde{\Sigma}$ divides the quotient of $\R^2 \times [-B, B]$ into two regions $\widetilde{R}_1$ and $\widetilde{R}_2$.
   \item Each of $\widetilde{R}_1$ and $\widetilde{R}_2$ contains at least one thread of $W$.
   \end{enumerate}
   We call such a surface $\widetilde{\Sigma}$ a \emph{separating plane} for the weave $W$.
   If there is no separating plane, a weave is said to be \emph{non-separable}.
\end{definition}

The following propositions show that a weave consists of topologically separated weave-connected components.

\begin{proposition}
  \label{proposition:3:weave-connected-component}
  If a weave consists of exactly one weave-connected component, then it is non-separable.
\end{proposition}

\begin{proposition}
  \label{proposition:3:connected-components-separable}
  If a weave consists of more than one weave-connected component, then it is separable.
\end{proposition}

We now prove Proposition~\ref{proposition:3:weave-connected-component}.
On the other hand, we will give a proof of Proposition~\ref{proposition:3:connected-components-separable} after we define the notion of a ``weave layer'' and discuss its properties in the next subsection.

The following lemma is used in the proof of Proposition~\ref{proposition:3:weave-connected-component}.
\begin{lemma}
  \label{lemma:3:alternating-quadruple}
  Let $W$ be a separable weave and $\widetilde{\Sigma}$ separating plane for $W$.
  Suppose that a subweave $(\{t^r_{i_1}, t^r_{i_2}\}, \{t^b_{j_1}, t^b_{j_2}\})_W$ of $W$ is divided by $\widetilde{\Sigma}$ into two parts.
  Then, the subweave does not form an \AltQuad.
\end{lemma}

\begin{proof}
  We assume that a weave $W$ is under Convention~\ref{convention:weave}.
  Let $W$ be a separable weave and $\widetilde{\Sigma}$ a separating plane for $W$.
  Then, there exists a positive number $B$ such that the threads of $W$ are contained in $\widetilde{U} = \R^2 \times [-B, B]$.
  We denote by $(U, \Sigma)$ the quotient of the pair $(\widetilde{U}, \widetilde{\Sigma})$ by periodic translations $\tau_1$ and $\tau_2$.
  By the definition of a separating plane, the surface $\Sigma$ divides $U$ into two regions $R_+$ and $R_-$ each of which is homeomorphic to $S^1 \times [0, 1]$.

  Let $\widetilde{\Sigma}_0$ be the $xy$-plane in $\R^3$, and $\Sigma_0$ the quotient of $\widetilde{\Sigma}_0$ by $\tau_1$ and $\tau_2$.
  So, $\Sigma_0$ separates $U$ into two regions $R^0_+$ and $R^0_-$, which are equal to $\Sigma_0 \times [0, B]$ and $\Sigma_0 \times [-B, 0]$, respectively.
  Hence, there exists an ambient isotopy, $f_t$, of $U$ that deforms $\Sigma$ into $\Sigma_0$, and (by interchanging the indices if necessary) $R_+$ and $R_-$ into $R^0_+$ and $R^0_-$, respectively.


  We denote by $\gamma^r_i$ and $\gamma^b_j$ the oriented simple closed curves in $U$ obtained by the quotient of the red thread $t^r_i$ and the blue thread $t^b_j$, respectively.
  (Remark that each thread admits a natural orientation as in Definition~\ref{definition:weave:sign}.)
  We have a projection $\pi\colon U \to \Sigma_0$ induced by the orthogonal projection $(x, y, z) \mapsto (x, y, 0)$.
  We can define the sign between the pair of $\gamma^r_i$ and $\gamma^b_j$ as similar to the sign $S(i, j)$ between $t^r_i$ and $t^b_j$, and it is consistent with $S(i, j)$.

  The links $f_t(\gamma^r_i)$ and $f_t(\gamma^b_j)$ project to the diagrams $\pi(f_t(\gamma^r_i))$ and $\pi(f_t(\gamma^b_j))$ on $\Sigma_0$, and any ambient isotopy between $f_t(\gamma^r_i)$ and $f_t(\gamma^b_j)$ in $U$ is realized by a finite sequence of Reidemeister moves and isotopies on $\Sigma_0$ (\cite[Lemma~1]{bourgoin}).
  So, the sum of the signs at the crossings between $\gamma^r_i$ and $\gamma^b_j$ is invariant under the ambient isotopy $f_t$.

  If $t^r_i$ and $t^b_j$ are separated by $\widetilde{\Sigma}$, then the sign $S(i, j)$ is determined by which of $f_1(\gamma^r_i)$ and $f_1(\gamma^b_j)$ is contained in $R^0_+$ and which is contained in $R^0_-$.
  In fact, if $f_1(\gamma^r_i)$ and $f_1(\gamma^b_j)$ are contained in $R^0_+$ and $R^0_-$, respectively, then $S(i, j)$ is $-1$.
  Otherwise, $S(i, j)$ is $+1$.


  We consider a subweave $(\{t^r_{i_1}, t^r_{i_2}\}, \{t^b_{j_1}, t^b_{j_2}\})_W$ of $W$ divided by $\widetilde{\Sigma}$ into two parts.
  Then, there are two cases: 
  \begin{enumerate}
    \item two corresponding curves of the same color are contained in either $R_+$ or $R_-$, and
    \item neither of them is on the same region.
  \end{enumerate}
  If there are two threads of the same color contained in either $R_+$ or $R_-$, then the sign between these two threads and a thread of the other color contained in the other region is the same.
  So, the subweave does not form an \AltQuad in this case.

  If there are no two threads of the same color contained in either $R_+$ or $R_-$, then each of the regions $R_+$ and $R_-$ contains one thread of each color.
  We can assume that $\gamma^r_{i_1}$ and $\gamma^b_{j_1}$ are contained in $R_+$, and $\gamma^r_{i_2}$ and $\gamma^b_{j_2}$ are contained in $R_-$.
  Then, $S(i_1, j_2) = -1 \neq +1 = S(i_2, j_1)$ holds.
  So, the subweave also does not form an \AltQuad in this case.
\end{proof}

\begin{proof}[Proof of Proposition~\ref{proposition:3:weave-connected-component}]
  We assume that a weave $W$ consists of exactly one weave-connected component.
  Let $\widetilde{\Sigma}$ be an embedded surface satisfying the conditions in Definition~\ref{definition:3:separable}, except for the last condition.
  By Lemma~\ref{lemma:3:alternating-quadruple}, $\widetilde{\Sigma}$ does not divide any \AltQuad into two parts.
  Since $W$ consists of exactly one weave-connected component, for any pair of threads $t$ and $t'$ in $W$, there exists a sequence of \AltQuads $U_0, \dots, U_l$ such that $t \in U_0$, $t' \in U_l$, and for each $k = 0, \dots, l-1$, $U_k$ and $U_{k+1}$ share a thread.
  So, the surface $\widetilde{\Sigma}$ does not divide any pair of threads in $W$ into two parts.
  Therefore, the weave $W$ is non-separable.
\end{proof}

\subsection{Weave layers}
\label{sec:weave:layers}

\par
The set of all weave-connected components admits a natural ordering based on their height.
Each weave-connected component is a non-separable sheet of interwoven threads and can be regarded as a sheet-like object.
Thus, for any two distinct components, one is intuitively positioned above the other.
The following definition formalizes this vertical relationship as the ``height order''.

\begin{figure}[htbp]
  \centering
  \begin{minipage}[b]{0.45\hsize}
    \centering
    \includegraphics[page=6, width=0.55\linewidth]{Figure/height_order_components.pdf}\\

    {\small
      $W_1 = (\{ t_1^r \}, \emptyset)_W$ or $(\emptyset, \{ t_1^b \})_W$\\
      $W_2 = (\{ t^r_2 \}, \emptyset)_W$ or $(\emptyset, \{ t^b_2 \})_W$
    }\\
    (a)
  \end{minipage}
  \begin{minipage}[b]{0.45\hsize}
    \centering
    \includegraphics[page=1, width=0.45\linewidth]{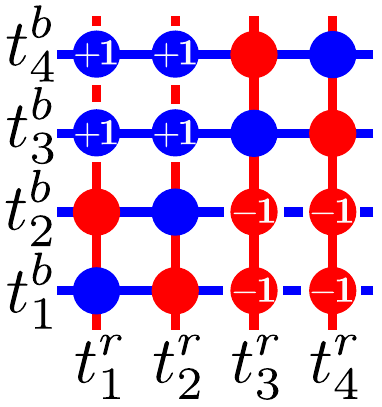}\\

    {\small
      $W_1 = (\{ t^{r}_1, t^{r}_2 \}, \{t^{b}_1, t^{b}_2 \})_W$\\
      $W_2 = (\{ t^{r}_3, t^{r}_4 \}, \{ t^{b}_3, t^{b}_4 \})_W$
    }\\
    (b) 
  \end{minipage}
\caption{Illustration of the height order $W_1 \preceq W_2$. (a) The case of parallel untangled threads. An intermediate thread determines the order.  (b) Non untangled thread case.}
  \label{fig:weaving:height_order}
\end{figure}

\begin{definition}
  \label{definition:weaving:height_order}
  Let $W_1$ and $W_2$ be two weave-connected components of a weave $W = (T^{r}, T^{b})$, and let $I^r_1, J^b_1$ and $I^r_2, J^b_2$ be the index sets of the red  and blue threads in $W_1$ and $W_2$, respectively.
  We denote by $S(i, j)$ the sign between the red thread $t^r_i$ and the blue thread $t^b_j$ in $W$.
  We say that $W_1 \preceq W_2$ if either one of the following conditions holds:
  \begin{enumerate}
    \item $W_1 = W_2$.
    \item[] The following conditions assume $W_1 \neq W_2$.
    \item 
      If $W_1$ and $W_2$ are untangled red threads, then there exists a blue thread $t^b_j \in W$ such that $S(i_1, j) = +1$ and $S(i_2, j) = -1$, where $I_1 = \{ i_1 \}$ and $I_2 = \{ i_2 \}$.
      Similarly, if $W_1$ and $W_2$ are untangled blue threads, then there exists a red thread $t^r_i \in W$ such that $S(i, j_1) = -1$ and $S(i, j_2) = +1$, where $J_1 = \{ j_1 \}$ and $J_2 = \{ j_2 \}$.
    \item 
      Suppose that $W_1$ or $W_2$ is not an untangled thread.
      Then, $S(i_1, j_2) = +1$ and $S(i_2, j_1) = -1$ hold for any $i_1 \in I_1$, $j_1 \in J_1$, $i_2 \in I_2$, and $j_2 \in J_2$.
  \end{enumerate}
  Illustrations of cases 2 and 3 are shown in Figure~\ref{fig:weaving:height_order}.
  As we will prove in Lemma~\ref{lemma:weaving:height order:partial order}, this relation defines a partial order on the set of weave-connected components.
  This partial order is referred to as the {\em height order}.
\end{definition}

\begin{remark}
  \label{remark:weaving:height order:either one}
  It is clear that if $W_1 \neq W_2$, and assuming that either $W_1$ or $W_2$ contains an \AltQuad, then $W_1 \preceq W_2$ and $W_2 \preceq W_1$ cannot both hold.
\end{remark}

\begin{remark}
  \label{remark:weaving:height order:total order}
  The height order is not a total order in general.
  This is because parallel untangled threads can exist in the same layer, making them incomparable.
\end{remark}

  An order relation between two components can always be determined, provided that at least one of them is not an untangled thread (see Lemma~\ref{lemma:weaving:height order:comparability}).
  Therefore, to establish a total order, we introduce the concept of a ``weave layer''.

\begin{figure}[htbp]
  \centering
  \includegraphics[page=3, width=0.30\linewidth]{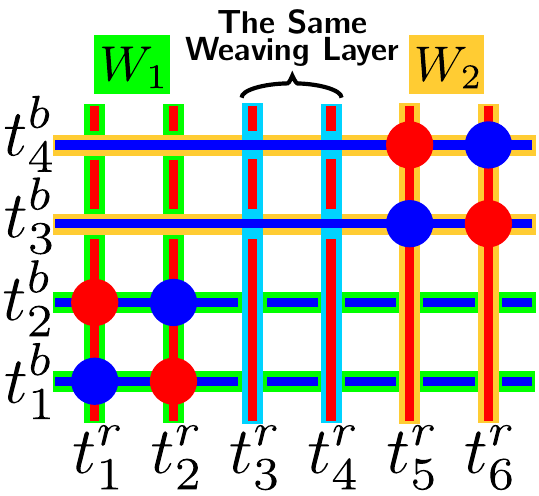}
  \caption{An example of weave layers. This weave has three weave layers: $W_1$ and $W_2$ are weave-connected components, while the remaining untangled threads form a single weave layer. The height order among these layers is $W_1 \preceq (\{ t^r_3, t^r_4 \}, \emptyset)_W \preceq W_2$.} 
  \label{fig:weaving:weaving_layer}
\end{figure}

\begin{definition}
  \label{definition:weaving:weaving_layer}
  We define a \emph{weave layer} as either
  \begin{enumerate}
    \item a weave-connected component that is not an untangled thread, or
    \item a maximal set of mutually incomparable untangled threads.
  \end{enumerate}
  An example of weave layers is illustrated in Figure~\ref{fig:weaving:weaving_layer}.
  The set consisting of these weave layers is then totally ordered by the height order.
\end{definition}

The properties of the height order for weave-connected components and weave layers can be summarized in the following two lemmas.

\begin{lemma}
  \label{lemma:weaving:height order:comparability}
  Any pair of weave-connected components are comparable under the height order,
  except in the case where both components are untangled threads with the same color.
\end{lemma}

\begin{lemma}
  \label{lemma:weaving:height order:partial order}
  The relation $\preceq$ defined in Definition \ref{definition:weaving:height_order} is a partial order on the set of weave-connected components of a weave.
  Furthermore, the relation $\preceq$ is a total order on the set of weave layers defined in Definition~\ref{definition:weaving:weaving_layer}.
\end{lemma}

According to the following strategy, we provide the proofs of Lemmas~\ref{lemma:weaving:height order:comparability} and~\ref{lemma:weaving:height order:partial order}.
First, we prove Lemma~\ref{lemma:weaving:height order:comparability} by using Lemmas~\ref{lemma:weaving:signs:single thread}, \ref{lemma:weaving:signs:two_parallel_threads}, and \ref{lemma:weaving:signs:two_intersecting_threads}.
We then prove Lemma~\ref{lemma:weaving:height order:partial order}.
However, since the proof of Lemma~\ref{lemma:weaving:height order:partial order} is straightforward but lengthy, we provide only the statement here and defer the proof to Appendix~\ref{sec:appendix:partial-order}.

\begin{lemma}
  \label{lemma:weaving:signs:single thread}
  Let $W_0$ be a weave-connected component of a weave $W$.
  If a thread $t$ of $W$ is not contained in $W_{0}$, then the signs between $t$ and the threads in $W_0$ are the same.
\end{lemma}

\begin{proof}
  Let $t$ be a thread in $W$ that is not in $W_0$.
  Without loss of generality, assume $t$ is a blue thread.

  \begin{claim}
    The two crossings between $t$ and the two red threads in an \AltQuad $U$ of $W_0$ have the same sign.
  \end{claim}
  \begin{claimproof}
    Let $U = (\{t^r_{i_1}, t^r_{i_2}\}, \{t^b_{j_1}, t^b_{j_2}\})_{W_0}$ be an \AltQuad of $W_0$.
    Since $U$ is a subweave of $W_0$ and $t$ is not in $W_0$, $t$ is distinct from both $t^b_{j_1}$ and $t^b_{j_2}$.
    We denote by $c_1$ and $c_2$ the crossings between $t$ and the red threads $t^{r}_{i_{1}}$ and $t^{r}_{i_{2}}$ in $U$, respectively.
    If the signs at $c_1$ and $c_2$ are different, then there exists a pair of crossings with opposite signs between $t$ and the red threads in $U$ by definition of \AltQuad (see Figure~\ref{fig:weaving:lemma:signs_component_thread}(a)).
    Hence, $t$ and $U$ form an \AltQuad, which contradicts the assumption that $t$ is not in $W_0$.
    So, the signs between $t$ and the two red threads in $U$ are the same.
  \end{claimproof}

  \begin{figure}[htbp]
    \begin{minipage}[b]{0.48\hsize}
      \centering
      \includegraphics[page=1, width=0.45\linewidth]{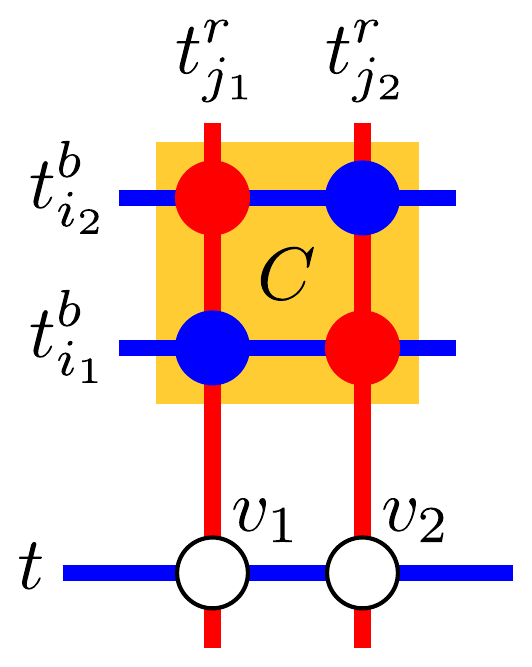}\\
      (a)
    \end{minipage}
    \begin{minipage}[b]{0.48\hsize}
      \centering
      \includegraphics[page=2, width=0.60\linewidth]{Figure/lemma_signs_component_thread.pdf}\\
      (b)
    \end{minipage}
    \caption{Illustration for the proof of the claims in Lemma~\ref{lemma:weaving:signs:single thread}. The figure~(a) shows the case where a thread $t$ intersects an \AltQuad $U$. The figure~(b) shows the case where two \AltQuads $U_1$ and $U_2$ share a blue thread.}
    \label{fig:weaving:lemma:signs_component_thread}
  \end{figure}

  \begin{claim}
    The signs between $t$ and the red threads in two intersecting \AltQuads $U_{1}$ and $U_{2}$ of $W_0$ are the same.
  \end{claim}

  \begin{claimproof}
    We first consider the case where $U_1$ and $U_2$ share a red thread. 
    By applying Claim~1 to $U_1$ and $U_2$, the signs between $t$ and the red threads in $U_1$ and $U_2$ are the same.
    This establishes the claim for this case.

    We next consider the case where $U_1$ and $U_2$ share a blue thread $t^b$. 
    We denote by $r_1$ and $r_1'$ the red threads in $U_1$, and by $r_2$ and $r_2'$ the red threads in $U_2$.
    Without loss of generality, we can assume that the signs between $t^b$ and the red threads $r_1$ and $r_2$ are signs $s_1$ and $s_2$, respectively, with $s_1 \neq s_2$ (see Figure~\ref{fig:weaving:lemma:signs_component_thread}(b)).
    Then, the signs between $t^b$ and the red threads $r_1'$ and $r_2'$ are $s_2$ and $s_1$, respectively.
    By Claim~1, the signs between $t$ and the red threads $r_1$ and $r_1'$ (resp. $r_2$ and $r'_2$) in $U_1$ (resp. $U_2$) are the same sign, which we denote by $s_3$ (resp. $s_4$). 

    We now proceed by contradiction.
    Assume that $s_3 \neq s_4$.
    By the definition of signs, $s_3$ must coincide with either $s_1$ or $s_2$.
    Without loss of generality, assume $s_3 = s_1$.
    Then, we have $s_4 \neq s_1$.
    Since $s_1 \neq s_2$, we have $s_4 = s_2$.
    Hence, a subweave $(\{r_1', r_2'\}, \{t^b, t\})_{W}$ forms an \AltQuad, which contradicts the assumption that $t \notin W_0$.
    Therefore, we must have $s_3 = s_4$.
\end{claimproof}

We take two red threads $r$ and $r'$ in $W_0$.
Since $W_0$ is weave-connected, there exists a sequence of \AltQuads $\{U_k\}_{k=0}^m$ such that $r \in U_0$ and $r' \in U_m$ and each successive pair $U_k$ and $U_{k+1}$ share a thread.
By repeatedly applying Claim~2 to $\{U_k\}_{k=0}^m$, the signs between the blue thread $t$ and the red threads $r$ and $r'$ are the same.
\end{proof}

\begin{lemma}
  \label{lemma:weaving:signs:two_parallel_threads}
  Let $t$ and $t'$ be both red threads or both blue threads in a weave-connected component $W_1$, and let $W_2$ be a weave-connected component that does not share any thread with $W_1$.
  Then, the signs between $t$ and the threads in $W_2$, and those between $t'$ and the threads in $W_2$ are the same.
\end{lemma}

\begin{proof}
  We only consider the case where $t$ and $t'$ are blue threads.
  By Lemma~\ref{lemma:weaving:signs:single thread}, the signs between $t$ and the red threads in $W_2$, and those between $t'$ and the red threads in $W_2$ are the same sign, respectively.
  On the other hand, since a red thread $t^r$ in $W_2$ is not in $W_1$, the sign between $t$ and $t^r$, and that between $t'$ and $t^r$ are the same sign.
  Therefore, the signs between $t$ and the red threads in $W_2$, and those between $t'$ and the red threads in $W_2$ are all identical.
\end{proof}

\begin{lemma}
  \label{lemma:weaving:signs:two_intersecting_threads}
  Let $W_{1}$ be a weave-connected component that contains an \AltQuad.
  Let $t^{r}$ and $t^{b}$ be a pair of a red thread and a blue thread that do not belong to $W_{1}$.
  Suppose that the signs between $t^{r}$ and the blue threads in $W_{1}$, and those between $t^{b}$ and the red threads in $W_{1}$ are the same.
  Then $t^{r}$ and $t^{b}$ are not lying in the same weave-connected component.
\end{lemma}

\begin{proof}
  Without loss of generality, we assume that the signs between $t^{r}$ and the blue threads $W_{1}$, and those between $t^{b}$ and the red threads $W_{1}$ are all $+1$.
  Assume that $t^{r}$ and $t^{b}$ are in the same weave-connected component $W_{2}$.

  We first show the sign between $t^{r}$ and $t^{b}$ is $+1$.
  Since $W_{1}$ contains an \AltQuad, there exists a pair of a red thread and a blue thread $(r_{1}, b_{1})$ in $W_{1}$ such that the sign between them is $-1$.
  Hence, if the sign between $t^{r}$ and $t^{b}$ is $-1$, then a subweave $(\{ t^{r}, r_1 \}, \{ t^{b}, b_1 \})_{W}$ forms an \AltQuad, which contradicts the assumption that both $t^{r}$ and $t^{b}$ are not in $W_{1}$.
  Hence, the sign between $t^{r}$ and $t^{b}$ is $+1$.

  The above discussion can be applied to any pair of a red thread and a blue thread in $W_{2}$.
  Thus, all signs between such pairs of threads are $+1$.
  This implies that $W_{2}$ does not contain any \AltQuad, which contradicts that $W_2$ is a weave-connected component containing two threads.
\end{proof}

\begin{proof}[Proof of Lemma~\ref{lemma:weaving:height order:comparability}]
  Let $W_1$ and $W_2$ be two weave-connected components.
  We assume that $W_1$ and $W_2$ have at least one crossing.
  This means that they are not both untangled threads without any crossings each other.

  If either $W_1$ or $W_2$ is an untangled thread, by Lemma~\ref{lemma:weaving:signs:single thread}, all crossings between the threads in $W_1$ and those in $W_2$ have the same sign.
  Thus, by Definition~\ref{definition:weaving:height_order}~(2), either $W_1 \preceq W_2$ or $W_2 \preceq W_1$ holds.

  We assume that both $W_1$ and $W_2$ contain \AltQuads.
  We can see by Lemmas~\ref{lemma:weaving:signs:single thread}, \ref{lemma:weaving:signs:two_parallel_threads} and \ref{lemma:weaving:signs:two_intersecting_threads} that $W_1$ and $W_2$ satisfy the condition in Definition~\ref{definition:weaving:height_order}~(3).
  Hence, either $W_1 \preceq W_2$ or $W_2 \preceq W_1$ holds.
\end{proof}

Next, we will show Proposition~\ref{proposition:3:connected-components-separable} by using the height order.

\begin{proof}
We assume that a weave $W$ has at least two weave-connected components.
If all components are untangled threads, then it is clear that $W$ is separable.
In the following, we consider the case where at least one component is not an untangled thread $W_0$.

We will divide the weave $W$ into two subweaves $W_-$ and $W_+$.
By Lemmas~\ref{lemma:weaving:height order:comparability} and~\ref{lemma:weaving:height order:partial order}, the set of weave-connected components of $W$ is partially ordered by the height order, and $W_0$ is comparable to any other component.
We denote by $W_-$ the subweave consisting of all components $W'$ such that $W' \preceq W_0$ and $W' \neq W_0$.
Furthermore, we denote by $W_+$ the subweave consisting of all components not contained in $W_-$.
We will show that there exists a separating plane between $W_-$ and $W_+$.
By the definition of the height order, we can move $W_+$ to a position with a large $z$-coordinate and sufficiently distant from $W_-$ using an ambient isotopy.
So, we can find a plane parallel to the $xy$-plane that separates $W_-$ and $W_+$.
The plane can be transformed by the inverse ambient isotopy to a separating plane between $W_-$ and $W_+$.
\end{proof}

We can see by using Proposition~\ref{proposition:3:connected-components-separable} that each distinct weave layer can also be separated by a (topological) plane.

\section{Stable configurations of weaves}
\label{sec:ode}
This section analyzes the stable configuration of a weave via the steepest descent method.
This technique is based on the principle of least action, which postulates that stable configurations in nature correspond 
to states where a relevant energy functional is minimized or at a local minimum.
By defining an appropriate functional to represent the physical phenomenon,
the method evolves the initial configuration along the negative gradient direction, 
which is the direction where the energy decreases most rapidly, 
to identify the configuration that attains an energy minimum.
A general discussion and illustrative examples of this method are provided in Appendix \ref{sec:appendix:a}.
\par
In the following, we first define the energy functional and discuss the fundamental properties of its associated steepest descent flow.
Subsequently, we perform an analysis for the case of weaves with a single layer (Section \ref{sec:eq:singlelayer}) 
and provide a proof of the main theorem (Theorem \ref{claim:weaving:ode:1-component}).
Finally, we present the analysis for the case of weaves with multiple layers (Section \ref{sec:weaving:2:2}) 
and establish the proof of the corresponding main theorem (Theorem \ref{claim:weaving:ode:untangle}).

\subsection{Energy functional for weaves}
\label{sec:eq:introduction}
In this subsection, we define the energy functional \eqref{eq:weaving:energy} 
and derive its associated steepest descent flow \eqref{eq:weaving:ode:00}.
First, we establish the energy inequality (Proposition \ref{claim:weaving:energy_bounds}),
which serves as the fundamental estimate for the analysis of the flow.
Furthermore, we provide several key results that play a crucial role in the subsequent sections,
including the commutativity of the Laplacians (Proposition \ref{claim:weaving:ode:matrix:1}) 
and the stationarity of the center of mass for each layer (Proposition \ref{claim:weaving:ode:center-eq}).
\par
Let $W=(T^r,T^b)$ be a weave.
A pair of non-parallel periodic translations $\tau_1$ and $\tau_2$ of $W$ induce linearly independent vectors $\iota_1, \iota_2 \in \mathbb{Z}\times\mathbb{Z}$, respectively.
Namely, for any crossing $c_{i,j}$ of $W$, we have $\tau_1(c_{i,j})=c_{(i,j) + \iota_1}$ and $\tau_2(c_{i,j})=c_{(i,j) + \iota_2}$.
Set $V = \{ (i, j) \in \Z \times \Z \mid (i, j) = \alpha \iota_1 + \beta \iota_2, 0 \leq \alpha < 1, 0 \leq \beta < 1 \}$.

By the definition of a weave, each red thread $t^r_i$ consists of vertices $v^r_{i,j}=(x_{i,j},y_{i,j},z^r_{i,j})$, and each blue thread $t^b_j$ consists of vertices $v^b_{i,j}=(x_{i,j},y_{i,j},z^b_{i,j})$.
These determine the following maps:
\begin{align*}
  \vfunc{x}&\colon \mathbb{Z}\times\mathbb{Z}\to\mathbb{R}^2, & 
  z_r&\colon \mathbb{Z}\times\mathbb{Z}\to\mathbb{R},&
  z_b&\colon \mathbb{Z}\times\mathbb{Z}\to\mathbb{R};\\
  \vfunc{x}&(i,j)=(x_{i,j},y_{i,j}), & z_r&(i,j)=z^r_{i,j}, & z_b&(i,j)=z^b_{i,j}.
\end{align*}
We call the triple $(\vfunc{x},z_r,z_b)$ a \emph{configuration} of the weave $W$.

\begin{definition}
  \label{definition:weaving:energy}
  For any index $u = (i, j) \in \Z \times \Z$, there exist unique real numbers $\alpha$, $\beta$ such that $u = \alpha \iota_1 + \beta \iota_2$.
  We define the projection $p\colon \Z \times \Z \to V$ and the shift $\vfunc{b}\colon \Z \times \Z \to \R^2$ by
  \begin{align*}
    p(i, j) &= (i, j) - \lfloor \alpha \rfloor \iota_1 - \lfloor \beta \rfloor \iota_2,\\
    \vfunc{b}(i, j) &= \lfloor \alpha \rfloor \tau_1 + \lfloor \beta \rfloor \tau_2.
  \end{align*}

  We define the energy $E(W)$ of the weave $W$ by
  \begin{equation}
    \label{eq:weaving:energy}
    \begin{aligned}
      E(W)
    &=
    \frac{1}{2} \sum_{v \in V} \sum_{u \in E_v} \lVert\vfunc{x}(p(u)) + \vfunc{b}(u) - \vfunc{x}(v)\rVert^2
    \\
    &+
    \frac{1}{2} \sum_{v \in V} \sum_{u \in E_v^r} \lvert z_r(p(u)) - z_r(v) \rvert^2
    \\
    &+
    \frac{1}{2} \sum_{v \in V} \sum_{u \in E_v^b} \lvert z_b(p(u)) - z_b(v) \rvert^2
    \\
    &+
    \sum_{v \in V} \frac{1}{\lvert z_b(v) - z_r(v) \rvert}.
    \end{aligned}
  \end{equation}
  Here, for $v:=(i,j)\in V$, we define subsets $E_v^r,\,E_v^b,\,E_v\subset V$ by $E_v^r=\{(i,j-1),(i,j+1)\}$, $E_v^b=\{(i-1,j),(i+1,j)\}$, and $E_v=E_v^r\cup E_v^b$.
  \par
  Throughout this section, we may also write $E(W)$ as $E(\Phi)$ using the configuration $\Phi=(\vfunc{x},z_r,z_b)$ of $W$.
\end{definition}
\par
\par
The energy functional \eqref{eq:weaving:energy} we define consists of two components: 
in the $xy$-direction, it is given by the sum of squared distances between adjacent vertices under the prescribed periodic conditions;
in the $z$-direction, it additionally incorporates repulsive interactions between corresponding vertices of different colors.
At a local minimum of this energy, 
the configuration achieves an equilibrium of forces in the $xy$-direction under the periodic constraints.
In the $z$-direction, this equilibrium is characterized by a balance between the attractive forces among adjacent vertices 
and the Coulomb repulsion between corresponding vertices of distinct colors.
\par
\begin{figure}[H]
  \centering
  \begin{tabular}{llll}
    (a)
    &(b)
    \\
    \includegraphics[bb=0 0 71 72,width=90pt]{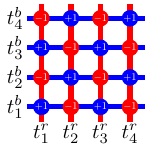}
    &\includegraphics[bb=0 0 71 72,width=90pt]{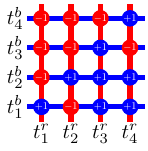}
    \\
    \mbox{}
    \\
    \includegraphics[bb=0 0 273 97,width=90pt]{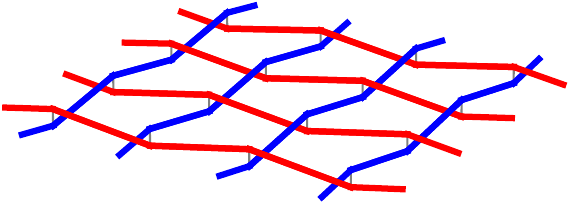}
    &\includegraphics[bb=0 0 272 202,width=90pt]{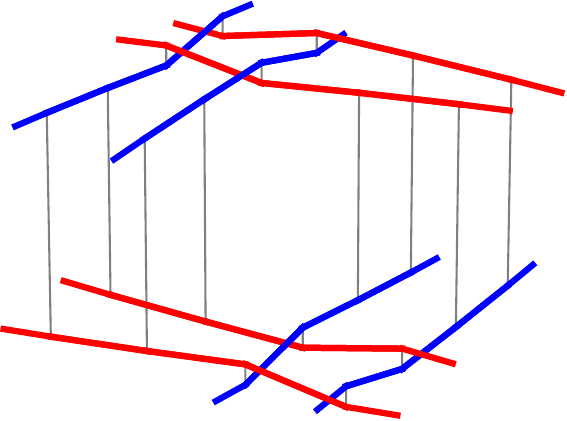}
    \\
  \end{tabular}
  \caption{
    (a) Single-layered configuration: 
    The solution to the steepest descent flow \eqref{eq:weaving:ode:00} 
    converges to the stable configuration shown on the right as $t\to\infty$.
    (b) Two-layered configuration: Each layer converges to a stable configuration, while the two layers diverge from one another.
  }
  \label{fig:weaving:0}
\end{figure}

We now define the notion of a stable configuration of a weave.
\begin{definition}
  \label{definition:ode:stable}
  A continuous deformation of configurations of a weave
  $\Phi(t)=(\vfunc{x}(t),z_{r}(t),z_{b}(t))$
  is said to be \emph{admissible} if it preserves the periodicity and the crossing signs of the orthogonal projection.
  \par
  Let $\Phi=(\vfunc{x},z_r,z_b)$ be a configuration of a weave $W$, and $C(\Phi)$ be 
  the set of configurations obtained from $\Phi$ by admissible deformations.
  We say that $\Phi_s \in C(\Phi)$ is a \emph{stable configuration} of $W$ if $\Phi_s$ is a local minimizer of the energy $E$ on $C(\Phi)$ (see Fig.~\ref{fig:weaving:0}).
\end{definition}
\par
To find the stable configuration of the weave, 
namely, the configuration that attains an energy minimum, 
we employ the steepest descent method, a standard approach in the analysis of variational problems.
Specifically, we consider the following differential equation \eqref{eq:weaving:ode:00},
which evolves the configuration $(\v{x}(t), \v{z}_r(t), \v{z}_b(t))$ in the direction of the steepest energy decrease.
For a given initial condition, we then analyze the asymptotic behavior of the solution to \eqref{eq:weaving:ode:00} as $t \to \infty$.
\par
The equation of the steepest descent flow for the energy (\ref{eq:weaving:energy}) is:
\begin{equation}
  \label{eq:weaving:ode:00}
  \begin{aligned}
    \frac{d}{dt} \vfunc{x}(t, v)
    &=
    \sum_{u \in E_v} (\vfunc{x}(t, p(u)) - \vfunc{x}(t, v)) + \sum_{u \in E_v} \vfunc{b}(u), 
    \\
    \frac{d}{dt} z_r(t, v)
    &=
    \sum_{u \in E_v^r}
    (z_r(t, p(u)) - z_r(t, v))
    - \frac{S(v)}{|z_b(t, v) - z_r(t, v)|^2}, 
    \\
    \frac{d}{dt} z_b(t, v)
    &=
    \sum_{u \in E_v^b}
    (z_b(t, p(u)) - z_b(t, v))
    + \frac{S(v)}{|z_b(t, v) - z_r(t, v)|^2}, 
    \\
  \end{aligned}
\end{equation}
where $S$ represents the signs at the crossings of the weave $W$.
\par
\par
We denote by $X$ the graph obtained from the orthogonal projection of a weave $W$ onto the $xy$-plane.
Namely, the vertex set of $X$ is the crossing points and the edge set consists of the images of the segments of the threads.
In the remainder of this section, we use the same symbols $X$, $T^r$, and $T^b$ to denote both the original graphs and their quotients by the translation symmetries.
\par
By using graph Laplacians, we can rewrite equation \eqref{eq:weaving:ode:00} in a simpler form.
Let $\Delta_X$, $\Delta_R$, and $\Delta_B$ be the graph Laplacians of $X$, $T^r$, and $T^b$, respectively.
Setting 
\begin{equation*}
  \bm{b} = \begin{bmatrix}
    \sum_{u \in E_v} \vfunc{b}(u)
  \end{bmatrix}_{v \in V} \text{ and }
  \v{r} = 
  \begin{bmatrix}
    \dfrac{S(v)}{|z_b(v) - z_r(v)|^2}
  \end{bmatrix}_{v \in V}, 
\end{equation*}
we may rewrite (\ref{eq:weaving:ode:00}) simply as
\begin{equation}
  \label{eq:weaving:ode:0}
  \begin{aligned}
    \frac{d}{dt} \v{x}
    &=
    \Delta_X \v{x} + \v{b}, 
    \\
    \frac{d}{dt} \v{z}_r
    &=
    \Delta_R \v{z}_r - \v{r}, 
    \\
    \frac{d}{dt} \v{z}_b
    &=
    \Delta_B \v{z}_b + \v{r}.
    \\
  \end{aligned}
\end{equation}
Here, similarly to $\bm{b}$ and $\bm{r}$, the symbols $\bm{x}$, $\bm{z}_r$, and $\bm{z}_b$ are vectors consisting of their corresponding functions.
\begin{remark}
  Since the first equation of \eqref{eq:weaving:ode:0} is linear,
  it is invariant under any linear transformation in the $xy$-direction. 
  Furthermore, as the second and third equations govern the $z$-direction, 
  they are also invariant under such transformations.
  Therefore, our results apply even if the basis vectors of the periodic lattice are not orthogonal.
\end{remark}
Here, we also remark that while the graph $X$ is quadrivalent, the graphs formed by $T^r$ and $T^b$ are divalent, 
namely, $T^r$ and $T^b$ are circular graphs.
Therefore, the Laplacian treated in this problem is based on the circular graph formed by the threads. 
The operators $\Delta_R$ and $\Delta_B$ appearing in the $\v{z}$-component of (\ref{eq:weaving:ode:0}) can be regarded as 
``one-dimensional Laplacians.'' 
On the other hand, since we deal with the collective deformation of layers as two-dimensional entities,
it is necessary to analytically recover the two-dimensional elements from these one-dimensional Laplacians.
\par
The energy inequality stated below is the most fundamental inequality in the analysis of steepest descent flows.
Under the assumption of the existence of a solution,
it provides a priori estimates that give the boundedness of terms such as $|\v{z}_r(t) - \v{z}_b(t)|$, 
and also ensures the convergence of the time-global solution.
\begin{proposition}
  \label{claim:weaving:energy_bounds}
  Let $\Phi(t)=(\vfunc{x}(t),z_r(t),z_b(t))$ be a time-dependent configuration of a weave.
  If $\Phi$ is a solution of (\ref{eq:weaving:ode:0}) on its maximal existence interval $J$, then for any $t \in J$,
  \begin{equation}
    \label{eq:weaving:energy_inequality}
    E(\Phi(t)) \le E(\Phi(0)) =: E_0.
  \end{equation}
\end{proposition}
A detailed proof of the energy inequality \eqref{eq:weaving:energy_inequality} in a generalized framework is presented
in Appendix \ref{sec:appendix:a}.
\par
The following proposition ensures, via energy inequalities, 
that the red and blue threads do not cross each other.
Consequently, it follows that $S$ is invariant along the flow.
The invariance of $S$ ensures that the over/under information
remains within the same smooth homotopy class throughout the evolution of the flow.
\begin{proposition}
  \label{claim:weaving:energy_bounds:1}
  Let $\Phi(t)=(\vfunc{x}(t),z_r(t),z_b(t))$ be a time-dependent configuration of a weave $W$.
  If $\Phi$ is a solution of (\ref{eq:weaving:ode:0}) on its maximal existence interval $J$, then for any $t \in J$,
  \begin{equation}
    \label{eq:weaving:e}
    |z_b(t, v) - z_r(t, v)| \ge E(\Phi)^{-1}
    \quad 
    \text{for all } v \in V.
  \end{equation}
  This implies that the signs at the crossings of $W$ are invariant along the solution, and that such deformation is admissible.
\end{proposition}
\begin{proof}
  By the definition of the energy, we have
  \begin{equation*}
    \frac{1}{|z_b(t, v) - z_r(t, v)|} \le  E(\Phi), \quad  \text{for all } v \in V,\, t \in J.
  \end{equation*}
  So, we obtain \eqref{eq:weaving:e}.
  Therefore, since the relative height of $z_b(t, v)$ and $z_r(t, v)$ does not change with time $t$, signs at each crossing are preserved.
\end{proof}
\par
\par
To find a stable configuration, we construct a solution to (\ref{eq:weaving:ode:0}) on $(0, \infty)$ 
for a given initial value $\Phi(0) = (\v{x}(0), \v{z}_r(0), \v{z}_b(0))$ 
and prove its convergence as $t \to \infty$. 
The first equation of (\ref{eq:weaving:ode:0}) is a linear heat equation 
with periodic boundary conditions arising from the periodicity of $W$ in the $xy$-direction.
Since it does not interact with $(\v{z}_r, \v{z}_b)$, it can be solved independently. 
As a result, in the $xy$-direction, the configuration converges to a state where the forces are balanced at each vertex.
This configuration is referred to as a harmonic realization (\cite{Kotani-Sunada}).
On the other hand, the second and third equations of (\ref{eq:weaving:ode:0}) 
interact with each other and contain the non-linear term $\v{r}$, 
which expresses the repulsive forces between the red vertices $(\vec{x}(v), z_r(v))$ 
and blue vertices $(\vec{x}(v), z_b(v))$.
\par
Since equation (\ref{eq:weaving:ode:0}) is Lipschitz continuous, 
the existence and uniqueness of a local-in-time solution are established for any given initial data (see \cite{Pontryagin}) 
\par
Propositions \ref{claim:weaving:ode:matrix:1} and \ref{claim:weaving:laplacian} play a fundamental role in the subsequent proofs.
\begin{proposition}
  \label{claim:weaving:ode:matrix:1}
  The Laplacians $\Delta_B$ and $\Delta_R$ are simultaneously diagonalizable. 
  The distinct eigenvalues of $\Delta_B$ and $\Delta_R$ are the same as those of the circular graphs $C_{N_b}$ and $C_{N_r}$, 
  respectively.
\end{proposition}
\begin{proof}
  Since taking a finite cover of a graph affects neither the simultaneous diagonalizability of the graph Laplacians nor their distinct eigenvalues, it suffices to consider weaves under Convention~\ref{convention:weave}.
  Hence, the Laplacians $\Delta_B$ and $\Delta_R$ are real symmetric matrices of order $n^2$.
  We index the vertices as shown in Figure \ref{fig:laplacian:1}
  \begin{figure}[H]
    \centering
    \includegraphics[bb=0 0 133 133,scale=0.5]{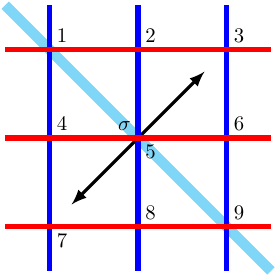}
    \quad
    \includegraphics[bb=0 0 218 218,scale=0.5]{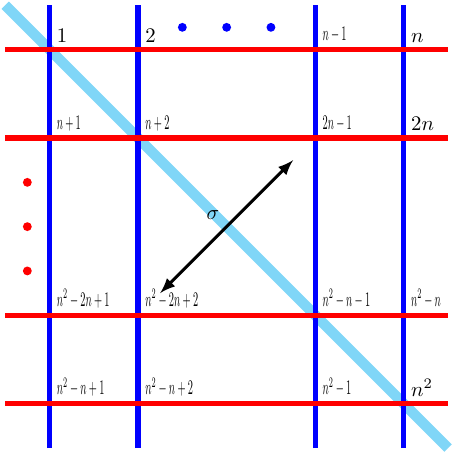}
    \caption{indexing and their exchange}
    \label{fig:laplacian:1}
  \end{figure}
  \par
  Let $\sigma \in \mathfrak{S}_{n^2}$ be a permutation defined by 
  $\sigma((k-1)n + j) = (j-1)n + k$ for $k, j = 1, \ldots, n$.
  Then $\Delta_R$ is obtained by permuting the rows and columns of $\Delta_B$ according to $\sigma$.
  Letting $T$ be the permutation matrix associated with $\sigma$, we obtain $\Delta_R = T \Delta_B T$.
  In the case of $n = 3$, $\sigma$ is 
  \begin{equation*}
    \begin{alignedat}{3}
      &\sigma(1) = 1, 
      &\quad
      &\sigma(2) = 4,
      &\quad
      &\sigma(3) = 7,
      \\
      &\sigma(4) = 2, 
      &\quad
      &\sigma(5) = 5, 
      &\quad
      &\sigma(6) = 8,
      \\
      &\sigma(7) = 3, 
      &\quad
      &\sigma(8) = 6, 
      &\quad
      &\sigma(9) = 9.
    \end{alignedat}
  \end{equation*}
  \par
  Let $\Delta_B = (b_1, \ldots, b_n)$ and $\Delta_R = (r_1, \ldots, r_n)$, where $b_i$, $r_j \in \R^n$.
  Since $\Delta_B$ and $\Delta_R$ are symmetric, the $(i, j)$-entries of their products are given by the inner products:
  \begin{equation*}
    \begin{aligned}
      (\Delta_B \Delta_R)_{ij} 
      &= \inner{b_i}{r_j}
      = \sum_{k=1}^{n^2} b_{ik} r_{jk}, 
      \\
      (\Delta_R \Delta_B)_{ij} 
      &= \inner{r_i}{b_j}
      = \sum_{k=1}^{n^2} r_{ik} b_{jk}, 
    \end{aligned}
  \end{equation*}
  since $\Delta_R$ and $\Delta_B$ are symmetric.
  \par
  The non-zero elements in the vector $b_i$ correspond to $b_{ii} = -2$ and the vertices adjacent to vertex $i$ along the blue thread. 
  The same applies to $r_j$. 
  Let $k$ be the index of the vertex where the blue thread containing $i$ intersects the red thread containing $j$. 
  Then we have $\langle b_i, r_j \rangle = b_{ik} r_{jk}$. 
  Similarly, let $\ell$ be the index of the vertex where the blue thread containing $j$ intersects the red thread containing $i$, so that $\langle r_i, b_j \rangle = r_{i\ell} b_{j\ell}$. 
  Since $\sigma(k) = \ell$, it follows that $r_{jk} = b_{j\ell}$ and $b_{ik} = r_{i\ell}$. 
  Consequently, $(\Delta_R \Delta_B)_{ij} = \langle r_i, b_j \rangle = r_{i\ell} b_{j\ell} = b_{ik} r_{jk} = \langle b_i, r_j \rangle = (\Delta_B \Delta_R)_{ij}$.
  Thus, $\Delta_B \Delta_R = \Delta_R \Delta_B$.
  \par
  The above argument holds for any indexing of vertices since any change in indexing can be represented by a conjugation with an appropriate permutation matrix.
\end{proof}
The following proposition provides the equations satisfied by the difference and the sum (average) of the height functions 
$\v{z}_r$ and $\v{z}_b$ along the flow, which will play a crucial role in the subsequent subsections.
\par
In the remainder, we simply write $d(t, v) = z_b(t, v) - z_r(t, v)$, $m(v) = z_b(t, v) + z_r(t, v)$, 
\begin{math}
  \v{d}(t) = 
  \begin{bmatrix}
    d(t, v)
  \end{bmatrix}_{v \in X}, 
\end{math}
and 
\begin{math}
  \v{m}(t) = 
  \begin{bmatrix}
    m(t, v)
  \end{bmatrix}_{v \in X}.
\end{math}
\begin{proposition}
  \label{claim:weaving:laplacian}
  For the Laplacians $\Delta_X$, $\Delta_B$, and $\Delta_R$, the following hold:
  \begin{align}
    \label{eq:weaving:laplacian:0}
    \Delta_B \v{z}_b + \Delta_R \v{z}_r
    &=
      \frac{1}{2} \Delta_X \v{m}
      + \frac{1}{2} (\Delta_B - \Delta_R)\v{d}, 
    \\
    \label{eq:weaving:laplacian:1}                                                                                                             
      \Delta_B \v{z}_b - \Delta_R \v{z}_r                                                                                                      
      &=                                                                                                                                       
      \frac{1}{2}(\Delta_B + \Delta_R)\v{d}                                                                                                    
      +                                                                                                                                        
      \frac{1}{2}(\Delta_B - \Delta_R)\v{m}.  
  \end{align}
\end{proposition}
\begin{proof}
  Using the commutativity $\Delta_B \Delta_R = \Delta_R \Delta_B$, we observe that 
  \begin{displaymath}
    \begin{aligned}
      \Delta_B \v{z}_b + \Delta_R \v{z}_r
      &=
      \frac{1}{2} (\Delta_B + \Delta_R)(\v{z}_b + \v{z}_r)
      \\
      &+ \frac{1}{2} \Delta_B \v{z}_b + \frac{1}{2} \Delta_R \v{z}_r
      - \frac{1}{2} \Delta_B \v{z}_r - \frac{1}{2} \Delta_R \v{z}_b
      \\
      &=
      \frac{1}{2} (\Delta_B + \Delta_R)(\v{z}_b + \v{z}_r)
      \\
      &+ \frac{1}{2} \Delta_B (\v{z}_b - \v{z}_r) 
      + \frac{1}{2} \Delta_R (\v{z}_r - \v{z}_b)
      \\
      &=
      \frac{1}{2} \Delta_X \v{m}                                                                                                               
      + \frac{1}{2} (\Delta_B - \Delta_R)\v{d}.  
    \end{aligned}
  \end{displaymath}
  Similarly, for the second identity:
  \begin{displaymath}
    \begin{aligned}
      \Delta_B \v{z}_b - \Delta_R \v{z}_r
      &=
      \frac{1}{2}(\Delta_B + \Delta_R)(\v{z}_b - \v{z}_r)
      \\
      &+ 
      \frac{1}{2} 
      \left(
        \Delta_B \v{z}_b
        - \Delta_R \v{z}_r
        - \Delta_R \v{z}_b
        + \Delta_B \v{z}_r
      \right)
      \\
      &=
      \frac{1}{2}(\Delta_B + \Delta_R)(\v{z}_b - \v{z}_r)
      +
      \frac{1}{2}(\Delta_B - \Delta_R)(\v{z}_b + \v{z}_r)
      \\                                                                                                                                       
      &=                                                                                                                                       
      \frac{1}{2}(\Delta_B + \Delta_R)\v{d}                                                                                                    
      +                                                                                                                                        
      \frac{1}{2}(\Delta_B - \Delta_R)\v{m}.   
    \end{aligned}
  \end{displaymath}
\end{proof}
In this paper, we analyze the relationship between layers by examining the barycenter of each layer.
\begin{definition}
  \label{definition:ode:center}
  Let $V$ be the index set of a weave $W$ defined above Definition~\ref{definition:weaving:energy}, and let $I^r$ and $J^b$ be the index sets of red and blue threads in $W$, respectively.
  For a subweave $W_1$ of $W$, we denote by $I^r_1 \subset I^r$ and $J^b_1 \subset J^b$ the index sets of red and blue threads in $W_1$, respectively.
  Let $V^r_1$ and $V^b_1$ be the index sets of the vertices in $W_1$ such that $V^r_1 = V \cap (I^r_1 \times J^b)$ and $V^b_1 = V \cap (I^r \times J^b_1)$.

  We define \emph{the barycenter (in the $z$-direction)} of $W_1$ as follows:
  \begin{displaymath}
    M_{W_1} 
    = \sum_{v \in V^r_1} z_r(v)
    + \sum_{v \in V^b_1} z_b(v),
  \end{displaymath}
  where $(\vfunc{x}, z_r, z_b)$ is a configuration of $W$.
\end{definition}
Next proposition establishes the invariance of the center of mass for each layer under the flow.
\begin{proposition}
  \label{claim:weaving:ode:center-eq}
  Let $\Phi(t)=(\vfunc{x}(t),z_r(t),z_b(t))$ be a time-dependent configuration of a weave $W$, and let $W_1$ be a subweave of $W$.
  Suppose that $\Phi$ is a solution of (\ref{eq:weaving:ode:0}).
  Then, we have
  \begin{equation}
    \label{eq:weaving:ode:center}
    \begin{aligned}
      &\frac{d}{dt} M_{W_1}(t) 
      = 
      \sum_{v \in V^b_1 \setminus V^r_1} 
      \frac{S(v)}{|z_b(v) - z_r(v)|^2}
      - 
      \sum_{v \in V^r_1 \setminus V^b_1} 
      \frac{S(v)}{|z_b(v) - z_r(v)|^2},
    \end{aligned}
  \end{equation}
  where $S$ represents the signs at the crossings of $W$.
  In particular, the barycenter $M_W(t)$ of $W$ does not depend on $t$.
\end{proposition}
\begin{proof}
  Throughout this proof, we adopt the notation introduced in Definition~\ref{definition:ode:center}.
  Summing (\ref{eq:weaving:ode:0}) over the red and blue threads in $W_1$, we obtain
  \begin{equation}
    \label{eq:weaving:ode:1:1}
    \begin{alignedat}{3}
      \frac{d}{dt} 
      \sum_{v \in V^b_1} z_b(v)
      &=
      \sum_{v \in V^b_1} \frac{S(v)}{|z_b(v) - z_r(v)|^2}, 
      \\
      \frac{d}{dt} 
      \sum_{v \in V^r_1} z_r(v)
      &=
      - 
      \sum_{v \in V^r_1} \frac{S(v)}{|z_b(v) - z_r(v)|^2}.
      \\
    \end{alignedat}
  \end{equation}
  Summing both equations in (\ref{eq:weaving:ode:1:1}), 
  we get 
  \begin{equation}
    \label{eq:weaving:ode:2:0}
    \begin{aligned}
      &\frac{d}{dt} 
      \sum_{v \in V^r_1 \cap V^b_1} (z_b(v) + z_r(v))
      + 
      \frac{d}{dt} 
      \sum_{v \in V^b_1 \setminus V^r_1} z_b(v)
      + 
      \frac{d}{dt} 
      \sum_{v \in V^r_1 \setminus V^b_1} z_r(v)
      \\
      &= \sum_{v \in V^b_1 \setminus V^r_1} 
      \frac{S(v)}{|z_b(v) - z_r(v)|^2}
      - 
      \sum_{v \in V^r_1 \setminus V^b_1} 
      \frac{S(v)}{|z_b(v) - z_r(v)|^2}.
    \end{aligned}
  \end{equation}
  If $W = W_1$, then 
  $V^b_1 \setminus V^r_1 = V^r_1 \setminus V^b_1 = \emptyset$, 
  hence the right hand side of (\ref{eq:weaving:ode:2:0}) is zero.
  Therefore the barycenter $M(t) = M_W(t)$ of $W$ is independent of $t$ along the solution of (\ref{eq:weaving:ode:00}).
\end{proof}
In the following, we assume $M(0) = 0$, which implies $M(t) = 0$ for all $t \in J$, without loss of generality.

\subsection{The case of weaves with a single layer}
\label{sec:eq:singlelayer}
In this subsection, we show that when the weave is single-layered, 
for any initial condition, there exists a time-global solution to the differential equation \eqref{eq:weaving:ode:00},
which converges to a stable configuration as $t \to \infty$ (Theorem \ref{claim:weaving:ode:1-component}).
As discussed in Section \ref{sec:eq:introduction},
the analysis can be reduced to the $z$-direction; specifically, it suffices to prove that 
$\v{z}_r(t)$ and $\v{z}_b(t)$ remain bounded for all $t>0$.
\par
First, by utilizing the single-layered structure, we establish the boundedness of $\v{d}(t) = \v{z}_b(t) - \v{z}_r(t)$.
In this step, the energy inequality (Proposition \ref{claim:weaving:energy_bounds}) 
and the properties of \AltQuads (Proposition \ref{claim:weaving:ode:minimal:d}) play an essential role.
\par
Furthermore, by combining the boundedness of $\v{d}(t)$ with the system of differential equations satisfied by 
$\v{d}(t)$ and $\v{m}(t) $(Proposition \ref{eq:weaving:ode:1:1}),
we demonstrate the boundedness of $\v{m}(t) = \v{z}_r(t) + \v{z}_b(t)$.
\par
In the proofs of the analytical inequalities presented below,
the symbol $C$ is used to denote generic positive constants.
It should be noted that these constants are not necessarily the same at each occurrence 
and may change from line to line.
Since we are investigating the existence of time-global solutions for given initial conditions,
it is essential that these constants $C$ are independent of the time variable $t$.
The following proposition shows that the difference $d$ in height functions remains bounded in the case of \AltQuad, 
playing an essential role in proving the convergence of a single-layered weave to a stable state.
\begin{proposition}
  \label{claim:weaving:ode:minimal:d}
  Let $W_1 = (\{ t^r_{i_1}, t^r_{i_2} \}, \{ t^b_{j_1}, t^b_{j_2} \})_W$ be an \AltQuad of a weave $W$.
  Then there exists a positive constant $C = C(E_0) > 0$, independent of $t \in J$, 
  such that 
  \begin{displaymath}
    \begin{aligned}
      &|d(t, (i_1, j_1))| + |d(t, (i_1, j_2))| + |d(t, (i_2, j_2))| + |d(t, (i_2, j_1))| \le C, 
      \\
      &\text{ where } d(t, v) = z_r(t, v) - z_b(t, v), 
    \end{aligned}
  \end{displaymath}
  along any solution of (\ref{eq:weaving:ode:00}).
\end{proposition}
\begin{proof}
  Let $W_1 = (\{ t^r_{i_1}, t^r_{i_2} \}, \{ t^b_{j_1}, t^b_{j_2} \})_W$ be an \AltQuad of a weave $W$.
  We denote by $v^r_{i, j}$ and $v^b_{i, j}$ the vertices of $t^r_i$ and $t^b_j$, respectively, defined as in Definition~\ref{definition:3:weave}.

  We take a closed polygonal path $c$ passing through the following vertices in order:
  $v^r_{i_1, j_1}$, 
  $v^b_{i_1, j_1}$, 
  $v^b_{i_2, j_1}$, 
  $v^r_{i_2, j_1}$, 
  $v^r_{i_2, j_2}$, 
  $v^b_{i_2, j_2}$, 
  $v^b_{i_1, j_2}$, 
  $v^r_{i_1, j_2}$, 
  and $v^r_{i_1, j_1}$. 
  The sum of the differences in the $z$-coordinates along the segments of $c$ is zero, namely,
  \begin{displaymath}
    \begin{aligned}
      0 
      &=
      (z_b(i_1, j_1) - z_r(i_1, j_1))
      + (z_b(i_2, j_1) -  z_b(i_1, j_1))
      \\
      &+ (z_r(i_2, j_1) - z_b(i_2, j_1))
      + (z_r(i_2, j_2) - z_r(i_2, j_1))
      \\
      &+ (z_b(i_2, j_2) - z_r(i_2, j_2))
      + (z_b(i_1, j_2) - z_b(i_2, j_2))
      \\
      &+ (z_r(i_1, j_2) - z_b(i_1, j_2))
      + (z_r(i_1, j_1) - z_r(i_1, j_2)).
    \end{aligned}
  \end{displaymath}
  Hence, we obtain
  \begin{equation}
    \label{eq:weaving:ode:1}
    \begin{aligned}
      &- d(i_1, j_1) + d(i_2, j_1) - d(i_2, j_2) + d(i_1, j_2)\\
      &=
      (z_b(i_2, j_1) -  z_b(i_1, j_1))
      + (z_r(i_2, j_2) - z_r(i_2, j_1))
      \\
      &+ (z_b(i_1, j_2) - z_b(i_2, j_2))
      + (z_r(i_1, j_1) - z_r(i_1, j_2)).
    \end{aligned}
  \end{equation}
  Taking the absolute value of both sides of (\ref{eq:weaving:ode:1}), we get
  \begin{equation}
    \label{eq:weaving:ode:2.1}
    \begin{aligned}
      &\lvert d(i_1, j_1) - d(i_2, j_1) + d(i_2, j_2) - d(i_1, j_2)\rvert\\
      &\le
      \lvert (z_b(i_2, j_1) -  z_b(i_1, j_1))\rvert
      + \lvert (z_r(i_2, j_2) - z_r(i_2, j_1))\rvert
      \\
      &+ \lvert (z_b(i_1, j_2) - z_b(i_2, j_2))\vert
      + \vert (z_r(i_1, j_1) - z_r(i_1, j_2))\rvert.
    \end{aligned}
  \end{equation}

  Each term on the right-hand side calculates the absolute value of the difference between two vertices that belong to the same thread.
  Hence, there exists a simple polygonal path along the thread connecting these two vertices.
  By the definition of the energy (\ref{eq:weaving:energy}), the absolute value of the difference along each segment is bounded by $E_0^{1/2}$.
  Therefore, there exists a positive integer $l$ such that the following inequality holds:
  \begin{equation}
    \label{eq:weaving:ode:2.2}
    \begin{aligned}
      &\lvert d(i_1, j_1) - d(i_2, j_1) + d(i_2, j_2) - d(i_1, j_2)\rvert
      \le
      l E_0^{1/2}.
    \end{aligned}
  \end{equation}

  Without loss of generality, we may assume that ($S(i_1, j_1)$, $S(i_2, j_1)$, $S(i_2, j_2)$, $S(i_1, j_2)$) = $(+1, -1, +1, -1)$, where $S$ represents the signs at the crossings of the weave $W$.
  Furthermore, for simplicity, we assume that $W$ satisfies the conditions of Convention~\ref{convention:weave}.
  Since $W_1$ is an \AltQuad, by Remark~\ref{remark:weave:sign}, we have
  \begin{equation*}
    \label{eq:weaving:ode:2.00}
    \begin{alignedat}{3}
      d(i_1, j_1) &= z_b(i_1, j_1) - z_r(i_1, j_1) > 0, &\quad d(i_2, j_1) &= z_r(i_2, j_1) - z_b(i_2, j_1) < 0, \\
      d(i_2, j_2) &= z_b(i_2, j_2) - z_r(i_2, j_2) > 0, &\quad d(i_1, j_2) &= z_r(i_1, j_2) - z_b(i_1, j_2) < 0.
    \end{alignedat}
  \end{equation*}
  Then, $d(i_1, j_1)$, $-d(i_2, j_1)$, $d(i_2, j_2)$, and $-d(i_1, j_2)$ have the same sign. So, we have
  \begin{equation}
    \label{eq:weaving:ode:20000}
    \begin{aligned}
      &\left|
        d(i_1, j_1) 
        - 
        d(i_2, j_1) 
        +
        d(i_2, j_2) 
        - 
        d(i_1, j_2) 
      \right|
      \\
      &=
      |d(i_1, j_1)| + |d(i_2, j_1)| + |d(i_2, j_2)| + |d(i_1, j_2)|. 
    \end{aligned}
  \end{equation}
  Therefore, by (\ref{eq:weaving:ode:2.2}) and (\ref{eq:weaving:ode:20000}), we obtain
  \begin{equation}
    \label{eq:weaving:ode:2}
    \begin{aligned}
      &|d(i_1, j_1)| + |d(i_2, j_1)| + |d(i_2, j_2)| + |d(i_1, j_2)|
      \le l E_0^{1/2}.
    \end{aligned}
  \end{equation}
\end{proof}
\begin{remark}
  Assume that a weave $W$ satisfies the conditions of Convention~\ref{convention:weave} and that a subweave $(\{ t^r_{i_1}, t^r_{i_2} \}, \{ t^b_{j_1}, t^b_{j_2} \})_W$ of $W$ is not an \AltQuad.
  For example, suppose that the signs ($S(i_1, j_1)$, $S(i_2, j_1)$, $S(i_2, j_2)$, $S(i_1, j_2)$) are given by $(+1, +1, -1, -1)$.
  Then, by Remark~\ref{remark:weave:sign}, we have
  \begin{equation}
    \label{eq:weaving:ode:2.001}
    \begin{alignedat}{3}
      & d(i_1, j_1) > 0, &\quad d(i_2, j_1) > 0, &\quad d(i_2, j_2) < 0, &\quad d(i_1, j_2) < 0.
    \end{alignedat}
  \end{equation}
  Then, we cannot establish (\ref{eq:weaving:ode:20000}), and thus we cannot obtain the uniform boundedness of $|\v{d}(t)|$.
  This assumption (\ref{eq:weaving:ode:2.001}) is an example of the untangled case (see Section \ref{sec:weaving:2:2}).
\end{remark}
By employing Proposition \ref{claim:weaving:ode:minimal:d},
one can show that $d$ remains bounded over the entire weave-connected component.
\begin{proposition}
  \label{claim:weaving:ode:maximal:d}
  Let $W_1$ be a weave-connected component of $W$.
  Then, there exists a positive constant $C = C(E_0, \{\lambda^X_k\}) > 0$, independent of $t \in J$, 
  such that 
  \begin{displaymath}
    \sum_{\substack{v \in V^r_1 \cap V^b_1}} |d(t, v)| \le C
  \end{displaymath}
  along any solution of (\ref{eq:weaving:ode:00}).
  Here, $V^r_1$ and $V^b_1$ are index sets defined in Definition~\ref{definition:ode:center}, and $\{ \lambda^X_k\}$ are eigenvalues of the graph Laplacian $\Delta_X$.
  In other words, the differences $d(t, v)$ for $v \in V^r_1 \cap V^b_1$ in a weave-connected component are uniformly bounded along the solution of (\ref{eq:weaving:ode:00}) in $t$.
\end{proposition}
\begin{proof}
  By the definition of weave-connectedness 
  there exists a sequence of \AltQuads $U_1, \ldots, U_k$ with 
  $U_k$ and $U_{k+1}$ share a thread (Definition \ref{definition:section2:altquad}).
  By Proposition \ref{claim:weaving:ode:minimal:d}, 
  $|d(t, v)|$ is bounded for each index of a thread pair $v$ in $U_k$.
  Since $U_k$ and $U_{k+1}$ share a thread, 
  $|d(t, v)|$ is also bounded for each index of a thread pair $v$ in $U_k \cup U_{k+1}$.
  Therefore, $|d(t, v)|$ is bounded for each index $v \in V^r_1 \cap V^b_1$.
\end{proof}
\par
The following proposition states a well-known property of ordinary differential equations;
however, we include it here for the reader's convenience and for the sake of completeness.
\begin{proposition}
  \label{claim:weaving:ode:solution}
  Consider the ordinary differential equation
  \begin{equation}
    \label{eq:weaving:ode:basic}
    \begin{aligned}
      \frac{d}{dt} \v{u}
      &=
      \Delta \v{u}
      + \v{f}
    \end{aligned}
  \end{equation}
  for a vector-valued function $\v{u}$.
  The solution of (\ref{eq:weaving:ode:basic}) is unique and is given by
  \begin{equation}
    \label{eq:weaving:ode:sol}
    \v{u}(t) = e^{t \Delta} \v{u}(0) + e^{t \Delta} \int_0^t e^{-s \Delta} \v{f}(s)\,ds
  \end{equation}
\end{proposition}
\begin{proof}
  Multiplying both sides of (\ref{eq:weaving:ode:basic}) by $e^{-t\Delta}$, 
  we obtain
  \begin{equation}
    \label{eq:weaving:ode:sol:1}
    \frac{d}{dt}\left(e^{-t \Delta} \v{u}(t)\right) = e^{-t \Delta} \v{f}(t).
  \end{equation}
  Integrating both sides of (\ref{eq:weaving:ode:sol:1}) over $[0, t]$ yields (\ref{eq:weaving:ode:sol}).
\end{proof}
Recall that $\Delta_X$, $\Delta_B$, and $\Delta_R$ are simultaneously  diagonalizable by Proposition~\ref{claim:weaving:ode:matrix:1}.
Let 
$0 = -\lambda^X_0 > -\lambda^X_1 \ge \cdots \ge -\lambda^X_{|V|-1}$, 
$0 = -\lambda^B_0 > -\lambda^B_1 \ge \cdots \ge -\lambda^B_{N_B-1}$, 
and 
$0 = -\lambda^R_0 > -\lambda^R_1 \ge \cdots \ge -\lambda^R_{N_R-1}$
be eigenvalues of $\Delta_X$, $\Delta_B$, and $\Delta_R$, respectively.
Let $\{\phi_k\}_{k=0}^{|V|-1}$ be an orthonormal basis of the vector space of functions on $V$ 
consisting of $L^2$-normalized eigenvectors of $\Delta_X$, 
and let $\{\phi_{k_j}\}$ be the set of eigenvectors corresponding to the non-zero eigenvalues of $\Delta_B - \Delta_X$.
\par
By employing these properties of eigenvalues and eigenfunctions, 
we can establish the uniform boundedness of $\v{m}(t)$ in the single-layered case.
In this proof, the uniform boundedness of $\v{d}(t)$ plays a crucial role.
\begin{proposition}
  \label{claim:weaving:ode:weaved:bounded}
  Assume
  a weave $W$ consists of a single weave layer.
  Then there exists a positive constant $C = C(E_0, \{\lambda^X_k\}) > 0$, independent of $t \in J$, 
  such that 
  $|\v{m}(t)| \le C$
  along any solution of (\ref{eq:weaving:ode:00}).
\end{proposition}
\begin{proof}
  From the assumption that $W$ consists of a single weave layer and (\ref{eq:weaving:laplacian:0}), 
  we obtain 
  \begin{equation}
    \label{eq:weaving:ode:weaved}
    \begin{aligned}
      \frac{d}{dt} \v{m}
      &=
      \frac{1}{2} \Delta_X \v{m}
      + \frac{1}{2} (\Delta_B - \Delta_R)\v{d}.
    \end{aligned}
  \end{equation}
  The expression of the solution of (\ref{eq:weaving:ode:weaved}) is easily obtained as
  \begin{equation}
    \label{eq:weaving:weaved:solution}
    \v{m}(t) = 
    e^{(t/2) \Delta_X} \v{m}(0)
    + 
    e^{(t/2) \Delta_X} 
    \int_0^t e^{-(s/2) \Delta_X} (\Delta_B - \Delta_R)\v{d}(s) \,ds.
  \end{equation}
  By using the eigenvector expansion, we may write 
  \begin{displaymath}
    \v{d}(t) = \sum_k \alpha_k(t) \phi_k
    = \sum_j \alpha_{k_j} (t) \phi_{k_j}, 
  \end{displaymath}
  and we obtain 
  \begin{equation}
    \label{eq:weaving:weaved:ode:estimate:0}
    \begin{aligned}
      &e^{(t/2) \Delta_X} 
      \int_0^t e^{-(s/2) \Delta_X} (\Delta_B - \Delta_R)\v{d}(s) \,ds
      \\
      &=
      e^
      {(t/2) \Delta_X} 
      \int_0^t \sum \lambda^{R-B}_{k_j}  \alpha_{k_j}(s) e^{-(s/2) \Delta_X}\phi_{k_j} \,ds
      \\
      &=
      \sum \lambda^{R-B}_{k_j} 
      \left(
        e^{-\lambda^X_{k_j}t/2} \int_0^t e^{\lambda^X_{k_j}s/2} \alpha_{k_j}(s)\,ds
      \right)
      \phi_{k_j}, 
    \end{aligned}
  \end{equation}
  where $\lambda^{R-B}_{k_j} = -\lambda^B_{k_j} + \lambda^R_{k_j} \not= 0$, 
  $\lambda^X_{k_j} = \lambda^B_{k_j} + \lambda^R_{k_j}$.
  Since we assume that the weave $W$ consists of a single layer, 
  by Proposition \ref{claim:weaving:ode:maximal:d}, 
  there exists a positive constant $C = C(E_0, \{\lambda^X_k\}) > 0$, independent of $t$, 
  such that 
  \begin{equation*}
    |\v{d}(t)|^2 = \sum |\alpha_k(t)|^2 \le C.
  \end{equation*}
  Hence we get 
  \begin{equation}
    \label{eq:weaving:weaved:ode:estimate:1}
    \begin{aligned}
      &
      \left|
        e^{(t/2) \Delta_X} 
        \int_0^t e^{-(s/2) \Delta_X} (\Delta_B - \Delta_R)\v{d}(s) \,ds
      \right|^2
      \\
      &\le
      \sum 
      C 
      \left|
        \lambda^{R-B}_{k_j} 
        e^{-\lambda^X_{k_j}t/2} \int_0^t e^{\lambda^X_{k_j}s/2}\,ds
      \right|^2.
    \end{aligned}
  \end{equation}
  Since eigenvalues of $\Delta_B$ and $\Delta_R$ are non-positive, 
  $\lambda^X_{k_j} = \lambda^B_{k_j} + \lambda^R_{k_j} = 0$ if and only if 
  $\lambda^B_{k_j} = \lambda^R_{k_j} = 0$, 
  hence $\lambda^{R-B}_{k_j} = 0$.
  Therefore, $\lambda^X_{k_j} \not= 0$ for each term in (\ref{eq:weaving:weaved:ode:estimate:1}).
  Moreover, we have 
  \begin{displaymath}
    \left|e^{-\lambda^X_{k_j}t/2} \int_0^t e^{\lambda^X_{k_j}s/2}\,ds\right|
    \le \left|\lambda^X_{k_j}\right|^{-1}, 
  \end{displaymath}
  yielding a uniform estimate 
  \begin{equation}
    \label{eq:weaving:weaved:ode:estimate:2}
    \begin{aligned}
      \left|
        e^{(t/2) \Delta_X} 
        \int_0^t e^{-(s/2) \Delta_X} (\Delta_B - \Delta_R)\v{d}(s) \,ds
      \right|^2
      \le C.
    \end{aligned}
  \end{equation}
  \par
  On the other hand, by Proposition \ref{claim:weaving:ode:center-eq} 
  and since the first eigenfunction $\phi_0$ of $\Delta_X$ is $\phi_0 = |V|^{-1/2} (1, \ldots, 1)$, 
  we have
  \begin{displaymath}
    \inner{\v{m}(t)}{\phi_0} 
    = 
    |V|^{-1/2} \sum_{v \in V} (z_b(t, v) + z_r(t, v))
    =
    |V|^{-1/2} M(t) = 0.
  \end{displaymath}
  Hence we obtain 
  \begin{equation}
    \label{eq:weaving:weaved:ode:estimate:3}
    \begin{aligned}
      \left|
        e^{(t/2) \Delta_X} \v{m}(t)
      \right| 
      \le C, 
    \end{aligned}
  \end{equation}
  and therefore we obtain a uniform estimate $|\v{m}(t)| \le C$ by 
  (\ref{eq:weaving:weaved:ode:estimate:0}), 
  (\ref{eq:weaving:weaved:ode:estimate:2}), and 
  (\ref{eq:weaving:weaved:ode:estimate:3}).
\end{proof}
Summarizing the above propositions, we obtain the main result for the single-layered case.
\begin{theorem}
  \label{claim:weaving:ode:1-component}
  If a weave $W$ consists of a single layer, 
  then the maximal interval of existence for the solution of (\ref{eq:weaving:ode:0}) 
  is $[0, \infty)$.
  Moreover, taking a suitable sequence $\{t_j\}$ with $t_j \to \infty$, 
  the solution $(\v{x}(t), \v{z}_b(t), \v{z}_r(t))$ converges to 
  a unique stationary solution $(\v{X}, \v{Z}_b, \v{Z}_r)$.
\end{theorem}
\begin{proof}
  By Propositions \ref{claim:weaving:ode:minimal:d} and \ref{claim:weaving:ode:weaved:bounded}, 
  the solution is contained in a compact set. 
  Therefore, the solution extends to $[0, \infty)$, 
  and by taking a suitable sequence, it converges to a stationary solution. 
  Since the functional $E$ is a sum of quadratic terms and terms of the form $x^{-1}$, 
  it is convex. 
  Hence, a stationary solution is unique in each homotopy class.
\end{proof}


\subsection{The case of weaves with multiple layers}
\label{sec:weaving:2:2}
In this subsection, we demonstrate that in the multi-layered case, the layers move apart at an order of $t^{1/3}$ 
while each individual layer remains bounded (Theorem \ref{claim:weaving:ode:untangle}).
In fact, unlike the single-layered case, the entire weave is not expected to converge to a stable structure 
since $\v{m}(t)$ remains bounded while $\v{d}(t)$ becomes unbounded (Proposition \ref{claim:weaving:ode:center_estimate}).
Nevertheless, regarding the $z$-coordinates of the centers of mass,
each layer remains within a bounded region (Theorem \ref{claim:weaving:ode:center_estimate:1}).
Furthermore, we show that the centers of mass of distinct layers diverge at an order of $t^{1/3}$
(Propositions \ref{claim:weaving:ode:lower_estimate} and \ref{claim:weaving:ode:upper_estimate}).
The following proposition shows that for the multiple-layered case, while $|\v{m}(t)|$ remains bounded, 
the divergence of $|\v{d}(t)|$ is estimated at order $t$; 
this result will be employed in the proofs of the subsequent propositions.
\begin{proposition}
  \label{claim:weaving:ode:center_estimate}
  Let $W$ be a weave consisting of multiple layers, and let $(\bm{x}, z_r, z_b)$ be a solution of (\ref{eq:weaving:ode:00}).
  Then there exists a positive constant $C = C(E_0, \{\lambda^X_q\}) > 0$, independent of $t \in J$, such that
  \begin{equation}
    \label{eq:weaving:ode:center_estimate}
    \begin{aligned}
      |\v{m}(t)| \le C, 
      \quad
      |\v{d}(t)| \le C(1+t), 
    \end{aligned}
  \end{equation}
  where $\{ \lambda^X_q \}$ are the eigenvalues of the graph Laplacian $\Delta_X$.
  In particular, $\v{z}_b(t)$ and $\v{z}_r(t)$ are bounded for all $t \in (0, T)$ and extend   
  to $(0, \infty)$.
\end{proposition}
\begin{proof}
  Even if $(T^r, T^b, S)$ is multiple layered, we have
  \begin{displaymath}
    \sup_{t \in J} |\v{m}(t)| \le C, 
  \end{displaymath}
  where the constant $C = C(E_0, \{\lambda^X_q\})$ is independent of $t \in J$.
  By (\ref{eq:weaving:laplacian:1}), $\v{d}$ and $\v{m}$ satisfy 
  \begin{displaymath}
    \frac{d}{dt}\v{d}
    = 
    \frac{1}{2} \Delta_X \v{d}
    +
    \frac{1}{2} (\Delta_B - \Delta_R) \v{m}
    + 
    \v{r}, 
  \end{displaymath}
  where 
  \begin{displaymath}
    \v{r}(t)
    =
    \begin{bmatrix}
      \dfrac{2 S(v)}{|z_b(t, v) - z_r(t, v)|^2}
    \end{bmatrix}_{v \in V}.
  \end{displaymath}
  Hence, $\v{d}$ is expressed as
  \begin{equation}
    \label{eq:weaving:untangle:proof:1:1}
    \begin{aligned}
      \v{d}(t)
      &=
      e^{(t/2) \Delta_X} \v{d}(0)
      +
      \frac{1}{2} e^{(t/2) \Delta_X} 
      \int_0^t e^{-(s/2) \Delta_X} (\Delta_B - \Delta_R) \v{m}(s)\, ds
      \\
      &+ 
      \frac{1}{2} e^{(t/2) \Delta_X} 
      \int_0^t e^{-(s/2) \Delta_X} \v{r}(s)\, ds.
    \end{aligned}
  \end{equation}
  To obtain the second estimate in the claim, 
  it is sufficient to show that $|\v{d}| \le C(t + 1)$ by using (\ref{eq:weaving:untangle:proof:1:1}).
  \par
  For the first term of the RHS of (\ref{eq:weaving:untangle:proof:1:1}), 
  since $\Delta_X$ is negative semi-definite, 
  the eigenvalues of $e^{(t/2)\Delta_X}$ are bounded from above by $1$.
  This ensures that the first term on the RHS of (\ref{eq:weaving:untangle:proof:1:1}) remains bounded.
  \par
  Similarly to the proof of Proposition \ref{claim:weaving:ode:weaved:bounded},
  and by using the commutativity of $\Delta_B$, $\Delta_R$, and $\Delta_X$, 
  the second term on the right-hand side of (\ref{eq:weaving:untangle:proof:1:1}) is also bounded. For the third term, we write 
  \begin{displaymath}
    \v{r}(t) = \sum_{k \ge 0} \alpha_k(t) \phi_k
  \end{displaymath}
  using an orthonormal system consisting of eigenvectors of $\Delta_X$.
  The $\phi_k$ components ($k \ge 1$) of the third term are bounded,
  and the $\phi_0$ component is estimated by 
  \begin{displaymath}
    \left|
      e^{(t/2) \Delta_X} 
      \int_0^t e^{-(s/2) \Delta_X} \v{r}(s)\, ds
    \right|
    \le C t.
  \end{displaymath}
  Therefore, we obtain 
  \begin{displaymath}
    \begin{aligned}
      |\v{m}(t)| \le C, \quad \text{and } 
      |\v{d}(t)| \le C(t + 1).
    \end{aligned}
  \end{displaymath}
\end{proof}
We consider a weave $W$ consisting of multiple weave layers $W_1, \ldots, W_K$.
By the definition of weave layers, each $W_k$ is either a single weave-connected component or consists of multiple untangled threads.
Let $W_{k, l}$ be a weave-connected component of $W_k$ for each $k = 1, \ldots, K$.
Let $V^r_{k, l}$ and $V^b_{k, l}$ be the index sets corresponding to $W_{k, l}$, as in Definition~\ref{definition:ode:center}.
Set
\begin{displaymath}
  w_{k, l} = \# (V^r_{k, l} \triangle V^b_{k, l}).
\end{displaymath}
Here, $X \triangle Y$ denotes the symmetric difference of sets $X$ and $Y$, i.e., $X \triangle Y = (X \setminus Y) \cup (Y \setminus X)$.
The following theorem shows that the height functions $\v{z}_r(t)$ and $\v{z}_b(t)$
remain bounded within each layer when viewed from its barycenter. 
This ensures that, in the subsequent proofs, it suffices to analyze the behavior of the barycenter of each layer.
\begin{theorem}
  \label{claim:weaving:ode:center_estimate:1}
  Assume that a weaving $W$ consists of weave layers $W_1, \ldots, W_K$ with a height order.
  Let $W_{k, l}$ be a weave-connected component of $W_k$ for each $k = 1, \ldots, K$.
  Let $V^r_{k, l}$ and $V^b_{k, l}$ be the index sets defined as Definition~\ref{definition:ode:center} for each weave-connected component $W_{k, l}$.
  Then, for any $v \in V^r_{k,l} \cup V^b_{k,l}$, 
  there exists a positive constant $C = C(E_0, \{\lambda^X_{q}\}, \{w_{k, l}\}) > 0$, 
  independent of $t \in J$, such that
  \begin{equation}
    \label{eq:ode:unweaving:difference}
    |z_b(t, v) - w_{k, l}^{-1} M_{W_{k,l}}(t)| \le C, 
    \quad
    |z_r(t, v) - w_{k, l}^{-1} M_{W_{k, l}}(t)| \le C
  \end{equation}
  along a solution of (\ref{eq:weaving:ode:00}). 
  In particular, if we take a subsequence $\{t_j\}$ with $t_j \to \infty$, 
  there exist $Z_b(v)$ and $Z_r(v)$ independent of $t$ such that 
  \begin{equation}
    \label{eq:ode:unweaving:convergence}
    \begin{aligned}
      &\lim_{t_j \to \infty} (z_b(t_j, v) - w_{k,l}^{-1} M_{W_{k,l}}(t_j)) = Z_b(v), 
      \\
      &\lim_{t_j \to \infty} (z_r(t_j, v) - w_{k,l}^{-1} M_{W_{k,l}}(t_j)) = Z_r(v).
    \end{aligned}
  \end{equation}
  Namely, $W_{k,l}$ converges under a shift by $M_{W_{k,l}}(t)$.
\end{theorem}
\begin{proof}
  We prove the case for $k=1$ and $W_1$ is a single weave-connected component.
  For simplicity, we denote by $W_1$ the weave-connected component $W_{1, 1}$.
  The other cases follow similarly.
  The index sets $V^r_1$ and $V^b_1$ are defined as in Definition~\ref{definition:ode:center} for $W_1$.
  We consider the equations for $\v{z}_b$, $\v{z}_r$, $\v{r}_b$, and $\v{r}_r$ with respect to vertices in $W_1$. 
  For example, let
  \begin{displaymath}
    \begin{bmatrix}
      \v{z}_b(t)
    \end{bmatrix}_{v \in V}
    =
    \begin{cases}
      z_b(t, v) &\text{ if } v \in V^b_1, \\
      0 &\text{ otherwise}.
    \end{cases}
  \end{displaymath}
  Then, the equation (\ref{eq:weaving:ode:0}) can be rewritten as 
  \begin{equation}
    \label{eq:weaving:ode:center:12.0}
    \begin{aligned}
      \frac{d}{dt} \v{z}_b
      &=
      \Delta_B \v{z}_b 
      + \v{r}_b, \\
      \frac{d}{dt} \v{z}_r
      &=
      \Delta_R \v{z}_r
      - \v{r}_r.
    \end{aligned}
  \end{equation}

  By the definition of the height order, the weave layer $W_1$ satisfies
  \begin{displaymath}
    \begin{aligned}
      S(v) = \begin{cases}
        +1 &\text{ if } v \in V^b_1 \setminus V^r_1, \\
        -1 &\text{ if } v \in V^r_1 \setminus V^b_1.
      \end{cases}
    \end{aligned}
  \end{displaymath}
  So, by Proposition~\ref{claim:weaving:ode:center-eq}, we obtain
  \begin{equation}
    \label{eq:weaving:ode:center:12}
    \begin{aligned}
      \frac{d}{dt} M_{W_1}(t) 
      &= 
      \sum_{v \in V^b_1 \setminus V^r_1} \frac{S(v)}{|z_b(v) - z_r(v)|^2}
      - 
      \sum_{v \in V^r_1 \setminus V^b_1} \frac{S(v)}{|z_b(v) - z_r(v)|^2}
      \\
      &=
      \sum_{v \in V^r_1 \triangle V^b_1} \frac{1}{|z_b(v) - z_r(v)|^2}.
    \end{aligned}
  \end{equation}
  For $v \in V^r_1 \cup V^b_1$, we define 
  \begin{displaymath}
    \begin{aligned}
      &\overline{z}_b(t, v) = z_b(t, v) - w_1^{-1} M_{W_1}(t), 
      \\
      &\overline{z}_r(t, v) = z_r(t, v) - w_1^{-1} M_{W_1}(t).
    \end{aligned}
  \end{displaymath}
  Then by (\ref{eq:weaving:ode:center:12.0}), (\ref{eq:weaving:ode:center:12}), 
  we obtain 
  \begin{equation}
    \begin{aligned}
      \label{eq:weaving:ode:center:13}
      \frac{d}{dt} \overline{\v{z}}
      &=
      \frac{d}{dt} \v{z}_b + \frac{d}{dt} \v{z}_r - 2 w^{-1}_1 \frac{d}{dt} M_{W_1}
      \\
      &=
      \Delta_B \overline{\v{z}}_b + \Delta_R \overline{\v{z}}_r  + \v{r}_b - \v{r}_r
      - 2 w^{-1}_1
      \sum_{v \in V^r_1 \triangle V^b_1} \frac{1}{|z_b(v) - z_r(v)|^2}.
    \end{aligned}
  \end{equation}
  Now we may rewrite linear terms in (\ref{eq:weaving:ode:center:13}) as 
  \begin{displaymath}
    \Delta_B \overline{\v{z}}_b + \Delta_R \overline{\v{z}}_r
    =
    \frac{1}{2} \Delta_X \overline{\v{m}} + \frac{1}{2}(\Delta_B - \Delta_R) \overline{\v{d}}. 
  \end{displaymath}
  Hence the equation (\ref{eq:weaving:ode:center:13}) is rewritten as
  \begin{equation}
    \label{eq:weaving:ode:center:14}
    \frac{d}{dt} \overline{\v{z}}
    =
    \frac{1}{2} \Delta_X \overline{\v{z}} + \frac{1}{2}(\Delta_B - \Delta_R) \overline{\v{d}}
    + \overline{\v{r}}, 
  \end{equation}
  where the nonlinear term $\overline{\v{r}}$ is given by 
  \begin{equation}
    \label{eq:weaving:ode:center:r}
    \overline{\v{r}}
    =
    \left\{
      \begin{aligned}
        &\displaystyle{
          \frac{1}{|z_b(v) - z_r(v)|^2} 
          - 2 w_1^{-1}
          \sum_{v \in V^r_1 \triangle V^b_1} \frac{1}{|z_b(v) - z_r(v)|^2}
        } 
        & \text{ for } v \in V^b_1 \setminus V^r_1, 
        \\
        &\displaystyle{
          -
          \frac{1}{|z_b(v) - z_r(v)|^2} 
          - 2 w_1^{-1}
          \sum_{v \in V^r_1 \triangle V^b_1} \frac{1}{|z_b(v) - z_r(v)|^2}
        }
        & \text{ for } v \in V^r_1 \setminus V^b_1,
        \\
        &0 & \text{otherwise}.
      \end{aligned}
    \right.
  \end{equation}
  We now prove the uniform boundedness of $\v{z}(t)$ using the expression for the solution to (\ref{eq:weaving:ode:center:14}):
  \begin{equation}
    \label{eq:weaving:ode:center:15}
    \begin{aligned}
      \overline{\v{z}}
      &=
      e^{(t/2)\Delta_X} \overline{\v{z}}(0)
      + 
      e^{(t/2)\Delta_X} 
      \int_0^t e^{-(s/2)\Delta_X} (\Delta_B - \Delta_R) \overline{\v{d}}(s)\,ds
      \\
      &+ 
      e^{(t/2)\Delta_X} 
      \int_0^t e^{-(s/2)\Delta_X}\overline{\v{r}}(s)\,ds.
    \end{aligned}
  \end{equation}
  Applying the same arguments used to establish the boundedness of (\ref{eq:weaving:untangle:proof:1:1}),
  we can conclude that the first term on the RHS of (\ref{eq:weaving:ode:center:15}) is bounded.
  \par
  We next show the boundedness of the second and third terms.
  From (\ref{eq:weaving:ode:center:r}), $\overline{\v{r}}$ is uniformly bounded.
  Since the zero-eigenspace of $\Delta_X$ is contained in the zero-eigenspace of $\Delta_B - \Delta_R$, 
  we have $\inner{\overline{\v{r}}(s)}{\phi_0} = 0$, where $\phi_0$ is the eigenvector of the zero-eigenvalue of $\Delta_X$.
  Therefore, the third term on the right-hand side of (\ref{eq:weaving:ode:center:15}) is bounded by similar calculations as in previous propositions.
  \par
  Finally, to prove the claim, we show the boundedness of $(\Delta_B - \Delta_R) \overline{\v{d}}(t)$.
  Since $W_1$ is a weave-connected component, 
  by the same arguments as in Proposition \ref{claim:weaving:ode:maximal:d},
  both $\overline{\v{d}}$ and $(\Delta_B - \Delta_R) \overline{\v{d}}(t)$ are uniformly bounded.
  Hence, $\overline{\v{z}}$ is uniformly bounded in $t$.
\end{proof}
In the following, we demonstrate that the center of mass of a given layer diverges 
from those of other layers at an order of $t^{1/3}$ (Proposition \ref{claim:weaving:ode:estimate}).
To this end, it is necessary to establish an estimate from below (Proposition \ref{claim:weaving:ode:lower_estimate}) 
and an estimate from above (Proposition \ref{claim:weaving:ode:upper_estimate}) separately.
\begin{proposition}
  \label{claim:weaving:ode:lower_estimate}
  Let $W$ be a weave consisting of multiple weave layers $W_1, \ldots, W_K$, and $W_{k,l}$ be a weave-connected component of $W_k$ for each $k = 1, \ldots, K$.
  Then there exist positive constants $C_0$ and $C_1$, independent of $t \in J$, such that
  \begin{equation}
    \label{eq:weaving:ode:lower_estimate}
    C_0 (t + C_1)^{1/3}
    \le
    |M_{W_{k,l}}(t) - M_{W_{k,l}^c}(t)|,
  \end{equation}
  along a solution of (\ref{eq:weaving:ode:00}).
  Here, $W_{k,l}^c$ is a subweave consisting of all threads in $W$ except for those in $W_{k,l}$.
\end{proposition}
\begin{proof}
  We prove the case for $k=1$ and $W_1$ is a single weave-connected component.
  The other cases follow similarly.
  We denote by $W_1$ the weave-connected component $W_{1, 1}$ for simplicity.
  \par
  Recall the equation
  \begin{equation}
    \label{eq:weaving:ode:1:6:1:11:1}
    \begin{aligned}
      \frac{d}{dt}(M_{W_1}(t) - M_{W_1^c}(t))
      &=
      \sum_{v \in V^r_1 \triangle V_1^b} \frac{2}{|z_b(v) - z_r(v)|^2}.
    \end{aligned}
  \end{equation}
  By Proposition \ref{claim:weaving:energy_bounds} and (\ref{eq:weaving:e}), $|z_b - z_r|^2 \le E_0$, so
  \begin{displaymath}
    |z_b(v) - z_r(v)|^2 \le E_0, 
    \quad
    |z_b(v) - z_r(v)|^{-2} \ge E_0^{-1}
  \end{displaymath}
  Therefore we obtain 
  \begin{displaymath}
    \begin{aligned}
      \frac{d}{dt}(M_{W_1}(t) - M_{W_1^c}(t))
      &\ge 
      2 w_1 E_0^{-1}, 
      \\
      M_{W_1}(t) - M_{W_1^c}(t)
      &\ge
      M_{W_1}(0) - M_{W_1^c}(0) 
      +
      2 w_1 E_0^{-1} t.
    \end{aligned}
  \end{displaymath}
  Hence there exists a $T > 0$ such that 
  if $t > T$, 
  \begin{displaymath}
    M_{W_1}(t) - M_{W_1^c}(t) > 0
  \end{displaymath}
  holds.
  In the following, we consider only on $t > T$.
  \par
  On the other hand, we have
  \begin{equation}
    \label{eq:weaving:lower_estimate:0}
    \begin{aligned}
      M_{W_1} - M_{W_1^c}
      &=
      \left(
        \sum_{v \in V^b_1} z_b(v)
        + 
        \sum_{v \in V^r_1} z_r(v)
      \right)
      - 
      \left(
        \sum_{v \in V \setminus V^b_1} z_b(v)
        + 
        \sum_{v \in V \setminus V^r_1} z_r(v)
      \right)
      \\
      &= 
      \sum_{v \in V^r_1 \cap V^b_1} (z_b(v) + z_r(v))
      -
      \sum_{v \in V \setminus (V^r_1 \cup V^b_1)} (z_b(v) + z_r(v))
      \\
      &+
      \sum_{v \in V^b_1 \setminus V^r_1} (z_b(v) - z_r(v))
      -
      \sum_{v \in V^r_1 \setminus V^b_1} (z_b(v) - z_r(v)).
    \end{aligned}
  \end{equation}
  Set 
  \begin{displaymath}
    M_{W_1 \triangle W^c_1}
    = 
    \sum_{v \in V^b_1 \setminus V^r_1} (z_b(v) - z_r(v))
    -
    \sum_{v \in V^r_1 \setminus V^b_1} (z_b(v) - z_r(v)).
  \end{displaymath}
  then by the definition of $W_1$, 
  we have
  \begin{displaymath}
    \begin{aligned}
      z_b(v) - z_r(v) > 0 &\quad\text{for } v \in V^b_1 \setminus V^r_1, \\
      z_b(v) - z_r(v) < 0 &\quad\text{for } v \in V^r_1 \setminus V^b_1, 
    \end{aligned}
  \end{displaymath}
  and 
  \begin{equation}
    \label{eq:weaving:lower_estimate:00}
    \begin{aligned}
      \frac{1}{M_{W_1 \triangle W_1^c}}
      &= 
      \frac{1}{
        \sum_{v \in V^b_1 \setminus V^r_1} (z_b(v) - z_r(v))
        -
        \sum_{v \in V^r_1 \setminus V^b_1} (z_b(v) - z_r(v))
      }
      \\
      &\le
      \sum_{v \in V^b_1 \setminus V^r_1}\frac{1}{|z_b(v) - z_r(v)|}
      + 
      \sum_{v \in V^r_1 \setminus V^b_1}\frac{1}{|z_b(v) - z_r(v)|}
      \\
      &=
      \sum_{v \in V^r_1 \triangle V_1^b}\frac{1}{|z_b(v) - z_r(v)|}.
    \end{aligned}
  \end{equation}
  Therefore, we obtain
  \begin{equation}
    \label{eq:weaving:ode:1:6:1:11:1:2}
    \begin{aligned}
      \frac{d}{dt}(M_{W_1}(t) - M_{W_1^c}(t))
      &=
      \sum_{v \in V^r_1 \triangle V_1^b} \frac{2}{|z_b(v) - z_r(v)|^2}
      \\
      &\ge
      \frac{1}{\lvert M_{W_1 \triangle W_1^c} \rvert}.
    \end{aligned}
  \end{equation}
  To show the claim, 
  we prove 
  \begin{equation}
    \label{eq:weaving:ode:16;1:1:11:1:2:1}
    \frac{1}{|M_{W_1 \triangle W_1^c}|}
    \ge
    \frac{1}{|M_{W_1} - M_{W_1^c}|}, 
    \text{ i.e., }
    |M_{W_1} - M_{W_1^c}|
    \ge
    |M_{W_1 \triangle W_1^c}|.
  \end{equation}
  \par
  Relationship between 
  $M_{W_1} - M_{W_1^c}$ and 
  $M_{W_1 \triangle W_1^c}$
  is 
  \begin{equation}
    \label{eq:weaving:ode:1:6:1:11:1:3}
    \begin{aligned}
      &M_{W_1}(t) - M_{W_1^c}(t)
      \\
      &= 
      \sum_{v \in V^r_1 \cap V^b_1} (z_b(v) + z_r(v))
      -
      \sum_{v \in V \setminus (V^r_1 \cup V^b_1)} (z_b(v) + z_r(v))
      + 
      M_{W_1 \triangle W_1^c}(t).
    \end{aligned}
  \end{equation}
  On the other hand, 
  by Theorem \ref{claim:weaving:ode:center_estimate:1}, 
  for any $v \in V^r_1 \cup V^b_1$, 
  there exist bounded functions $C_{W_1, b}(v)$ and  $C_{W_1, r}(v)$ 
  such that 
  \begin{displaymath}
    \begin{aligned}
      z_b(t, v) &= C_{W_1, b}(v) + w_1^{-1} M_{W_1}(t), & z_r(t, v) &= C_{W_1, r}(v) + w_1^{-1} M_{W_1}(t).
    \end{aligned}
  \end{displaymath}
  Hence, there exists positive constant $c_{W_1} = c_{W_1}(E_0, \{\lambda^X_q\}, \{w_k\}) > 0$ and 
  $c_{W_1^c} = c_{W_1^c}(E_0, \{\lambda^X_q\}, \{w_k\}) > 0$ such that
  \begin{displaymath}
    \begin{aligned}
      \sum_{v \in V^r_1 \cap V^b_1} (z_b(v) + z_r(v))
      &= 
      c_{W_1} M_{W_1}(t) + \sum_{v \in V^r_1 \cap V^b_1} (C_{W_1, b}(v) + C_{W_1, r}(v)), 
      \\
      \sum_{v \in V \setminus (V^r_1 \cup V^b_1)} (z_b(v) + z_r(v))
      &= 
      c_{W_1^c} M_{W_1^c}(t) + \sum_{v \in V\setminus (V^r_1 \cup V^b_1)} (C_{W_1^c, b}(v) + C_{W_1^c, r}(v)), 
    \end{aligned}
  \end{displaymath}
  and we obtain 
  \begin{equation}
    \label{eq:weaving:lower_estimate:1}
    \begin{aligned}
      &\sum_{v \in V^r_1 \cap V^b_1} (z_b(v) + z_r(v))
      -
      \sum_{v \in V \setminus (V^r_1 \cup V^b_1)} (z_b(v) + z_r(v))
      \\
      &= 
      c_{W_1} M_{W_1}(t) 
      -
      c_{W_1^c} M_{W_1^c}(t) 
      +
      \sum_{v \in V^r_1 \cap V^b_1} (C_{W_1, b}(v) + C_{W_1, r}(v)), \\
      &- \sum_{v \in V \setminus (V^r_1 \cup V^b_1)} (C_{W_1^c, b}(v) + C_{W_1^c, r}(v)).
    \end{aligned}
  \end{equation}
  Hence there exist positive constants $C_0 = C_0(E_0, \{\lambda^X_q\}, \{w_k\})> 0$ 
  and $C_1 = C_1(E_0, \{\lambda^X_q\}, \{w_k\})> 0$, which are independent of $t \in J$, 
  such that
  \begin{equation}
    \label{eq:weaving:lower_estimate:2}
    \begin{aligned}
      &
      \left|
        \sum_{v \in V^r_1 \cap V^b_1} (z_b(v) + z_r(v))
        -
        \sum_{v \in V \setminus (V^r_1 \cup V^b_1)} (z_b(v) + z_r(v))
      \right|
      \\
      &\le
      C_0|M_{W_1}(t) - M_{W_1^c}(t)|
      \\
      &+ 
      \left|
        \sum_{v \in V^r_1 \cap V^b_1} (C_{W_1, b}(v) + C_{W_1, r}(v)) 
        - \sum_{v \in V \setminus (V^r_1 \cup V^b_1)} (C_{W_1^c, b}(v) + C_{W_1^c, r}(v))
      \right|
      \\
      &\le
      C_0|M_{W_1}(t) - M_{W_1^c}(t)| + C_1.
    \end{aligned}
  \end{equation}
  Therefore by (\ref{eq:weaving:ode:1:6:1:11:1:3}) and (\ref{eq:weaving:lower_estimate:1}), 
  we obtain 
  \begin{equation}
    \label{eq:weaving:ode:1:6:1:11:1:4}
    \begin{aligned}
      \left|
      M_{W_1 \triangle W_1^c}
      \right|
      &\le
      \left|
        M_{W_1} - M_{W_1^c} 
      \right|
      +   
      C_0|M_{W_1}(t) - M_{W_1^c}(t)| + C_1
      \\
      &=
      C_0 \left|
        M_{W_1} - M_{W_1^c} 
      \right|
      + C_1.
    \end{aligned}
  \end{equation}
  Now combining (\ref{eq:weaving:ode:1:6:1:11:1:2}) and (\ref{eq:weaving:ode:1:6:1:11:1:4}), 
  we obtain 
  \begin{equation}
    \label{eq:weaving:ode:1:6:1:11:1:5}
    \begin{aligned}
      \frac{d}{dt}(M_{W_1}(t) - M_{W_1^c}(t))
      &\ge
      \frac{1}{|M_{W_1 \triangle W_1^c}|}
      \ge
      \frac{C_0}{|M_{W_1} - M_{W_1^c}| + C_1}
    \end{aligned}
  \end{equation}
  for positive constants $C_i = C_i(E_0, \{\lambda^X_q\}, \{w_k\}) > 0$, $i = 1,\, 2$.
  Therefore we obtain (\ref{eq:weaving:ode:lower_estimate}).
\end{proof}
\begin{proposition}
  \label{claim:weaving:ode:upper_estimate}
  Let $W$ be a weave consisting of multiple weave layers $W_1, \ldots, W_K$ with height order, and $W_{k,l}$ be a weave-connected component of $W_k$ for each $k = 1, \ldots, K$.
  Then there exist positive constants $C_0 = C_0(E_0, \{\lambda^X_{q}\}, \{w_{k,l}\}) > 0$ 
  and $C_1 = C_1(E_0, \{\lambda^X_{q}\}, \{w_{k,l}\}) > 0$, 
  independent of $t \in J$, 
  such that 
  \begin{equation}
    \label{eq:weaving:ode:upper_estimate}
    |M_{W_{k, l}}(t) - M_{W_{k, l}^c}(t)| 
    \le
    C_0 (t + C_1)^{1/3}, 
  \end{equation}
  holds along a solution of (\ref{eq:weaving:ode:00}).
\end{proposition}
\begin{proof}
  We provide the proof for $k=1$ and $W_1$ is a single weave-connected component.
  The cases for other $k$ follow similarly.
  We denote by $W_1$ the weave-connected component $W_{1, 1}$ for simplicity.
  \par
  Recall that the governing equation is
  \begin{equation}
    \label{eq:weaving:ode:1:6:1:11:1:1}
    \begin{aligned}
      \frac{d}{dt}(M_{W_1}(t) - M_{W_1^c}(t))
      &=
      \sum_{v \in V^r_1 \triangle V_1^b} \frac{2}{|z_b(v) - z_r(v)|^2}. 
    \end{aligned}
  \end{equation}
  By Proposition \ref{claim:weaving:ode:lower_estimate}, 
  we have
  \begin{displaymath}
    c_0 (t + c_1)^{1/3}
    \le
    |M_{W_1}(t) - M_{W_1^c}(t)|. 
  \end{displaymath}
  Furthermore, we observe that
  \begin{displaymath}
    \begin{aligned}
      c_0(t + c_1)^{1/3} 
      &\le 
      w_1^{1/3} |M_{W_1}(t) - M_{W_1^c}(t)|
      \\
      &\le
      |z_b(t) - z_r(t)
      +
      |z_b(t) - w_1^{-1} M_{W_1}(t)|
      +
      |z_r(t) - w_1^{-1} M_{W_1}(t)|.
    \end{aligned}
  \end{displaymath}
  On the other hand, 
  Theorem \ref{claim:weaving:ode:center_estimate:1} implies that
  \begin{displaymath}
    \begin{aligned}
      &|z_b(t, v) - w_1^{-1} M_{W_1}(t)| \le C, 
      \\
      &|z_r(t, v) - w_1^{-1} M_{W_1}(t)| \le C. 
    \end{aligned}
  \end{displaymath}
  Thus, for sufficiently large $t$, 
  we obtain 
  \begin{displaymath}
    \begin{aligned}
      C_0(t + C_1)^{1/3} 
      &\le w_1^{-1} |M_{W_1}(t) - M_{W_1^c}(t)|
      \le |z_b(t, v) - z_r(t, v)|.
    \end{aligned}
  \end{displaymath}
  Therefore, it follows that
  \begin{displaymath}
    \frac{1}{|z_b(t, v) - z_r(t, v)|^2} \le C_0(t + C_1)^{-2/3}
  \end{displaymath}
  and 
  \begin{equation}
    \label{eq:weaving:ode:1:6:1:11:1:20}
    \begin{aligned}
      \frac{d}{dt}(M_{W_1}(t) - M_{W_1^c}(t))
      &=
      \sum_{v \in V^r_1 \triangle V_1^b} \frac{2}{|z_b(v) - z_r(v)|^2}
      \\
      &\le
      C_0(t + C_1)^{-2/3}.
    \end{aligned}
  \end{equation}
  Integrating (\ref{eq:weaving:ode:1:6:1:11:1:20}) yields
  \begin{displaymath}
    |M_{W_1}(t) - M_{W_1^c}(t)|
    \le C_0(t + C_1)^{1/3}.
  \end{displaymath}
\end{proof}
\begin{proposition}
  \label{claim:weaving:ode:estimate}
  Let $W$ be a weave consisting of multiple weave layers $W_1, \ldots, W_K$ with height order, and $W_{k,l}$ be a weave-connected component of $W_k$ for each $k = 1, \ldots, K$.
  Then the mean values satisfy
  \begin{equation}
    |M_{W_{k, l}}(t) - M_{W_{k, l}^c}(t)| 
    \sim O(t^{1/3}), 
    \quad
    \text{ as } 
    t \to \infty, 
  \end{equation}
  along a solution of (\ref{eq:weaving:ode:00}).
\end{proposition}
Summarizing the above propositions, we obtain the main result for the multi-layered case.
\begin{theorem}
  \label{claim:weaving:ode:untangle}
  Let $W$ be a weave consisting of multiple weave layers $W_1, \ldots, W_K$ with height order, and $W_{k,l}$ be a weave-connected component of $W_k$ for each $k = 1, \ldots, K$.
  Then the weave-connected component $W_{k,l}$ moves apart at a rate of order $t^{1/3}$ along the steepest descent flow (\ref{eq:weaving:ode:0}) 
  of the energy (\ref{eq:weaving:energy}). 
  Namely, if we take a weave-connected component $W_{k_1,l}$ and a weave layer $W_{k_2}$ with $k_1 \neq k_2$,
  then the distance between their mean values $M_{W_{k_1, l}}(t)$ and $M_{W_{k_2}}(t)$ increases at order $t^{1/3}$, 
  while each weave-connected component $W_{k, l}(t)$ converges to a neighborhood of the plane $z = M_{W_{k, l}}(t)$ at the same order.
\end{theorem}


\appendix
\section{Appendix: Ideas and methods of variational problems}
\label{sec:appendix:a}
In this appendix, we summarize the fundamental ideas and methods of variational problems used in the main text. 
The calculus of variations is a standard framework for identifying stable configurations 
when an abstract object is realized in Euclidean space. 
In natural phenomena,
the principle of least action dictates that a system selects a configuration that minimizes,
or at least renders stationary, a certain energy functional. 
Such configurations are stable with respect to local deformations.
For precise terminology and theorems, see \cite{Jost}.
\par
Given a system (such as graphs or manifolds),
we define an energy functional $E \colon X \to \R$ on a functional space $X$.
The state realized by the system is typically given by a solution with minimal energy.
Such a solution is called a locally energy-minimizing solution.
This approach is known as the variational principle. 
In the case where the vector space $X$ is finite-dimensional,
a solution with minimal energy exists if the energy functional $E$ satisfies the following conditions:
a) $E$ is bounded from below, and 
b) any sub-level set $L(c) = \{x \in X \mid E(x) \le c\}$ is compact.
Furthermore, if $E$ is a convex function on the entire space $X$,
then the minimal energy solution is unique and is the global energy minimizer.
\par
The {\em steepest descent method} is widely known as a technique for finding energy-minimizing solutions.
This method employs evolutionary differential equations where the energy decreases monotonically along the solution trajectory.
When conditions a) and b) above hold,
there exists a smooth solution for all time $t \in [0, \infty)$ for any initial condition, 
which converges to a stationary solution with minimal energy.
\par
On the other hand, when an energy functional is bounded from below, the following two questions are of importance:
A) Does the evolution equation admit a unique time-global solution for any initial value?
B) If so, does the solution converge to a stationary solution?
If both A and B are answered affirmatively, we can demonstrate the existence of a locally energy-minimizing solution.
That is to say, the solution of the steepest descent flow,
starting from an initial condition with a small perturbation of a stationary solution,
converges back to the stationary solution. In this sense, a minimal energy solution is stable under small perturbations.
\par
More precisely, let $E \colon X \to \R$ be a smooth functional on a functional space $X$.
For a family of perturbations $\{x_s\}_{s \in (-\epsilon, \epsilon)} \subset X$ around $x_0 \in X$, we have 
\begin{equation}
  \label{eq:a:00}
  \left.\frac{d}{ds} E(x_s)\right|_{s=0}
  = 
  \inner{\nabla E(x_0)}{v}, 
  \quad
  \text{where }
  v = \left.\frac{d x_s}{ds}\right|_{s=0}.
\end{equation}
Therefore, the steepest descent flow of the energy $E$ for $x(t) \in X$, $t \in [0, T)$, given by
\begin{equation}
  \label{eq:a:01}
  \frac{d}{dt}x(t)
  =
  -\nabla E(x(t))
\end{equation}
is the most efficient way to reduce energy. 
In fact, taking the inner product with $x'(t)$ on both sides of (\ref{eq:a:01}), we obtain 
\begin{equation}
  \label{eq:a:02}
  |x'(t)|^2
  = 
  -\inner{\nabla E(x(t))}{x'(t)}.
\end{equation}
By using (\ref{eq:a:00}) and (\ref{eq:a:02}) and integrating over $[0, T]$, we obtain 
\begin{equation}
  \label{eq:a:03}
  \int_0^T |x'(t)|^2\,dt
  =
  -\int_0^T \frac{d}{dt} E(x(t))\,dt
  =
  E(x(0)) - E(x(T)).
\end{equation}
for any $T > 0$.
Consequently, as long as a solution $x(t)$ of (\ref{eq:a:01}) exists,
we have the {\em energy inequality}:
\begin{equation}
  \label{eq:a:04}
  E(x(t)) \le E(x(0)), \quad \text{for } t > 0.
\end{equation}
\par
Now, assume the energy $E$ is bounded from below and any sub-level set $L(c) = \{x \in X \mid E(x) \le c\}$ is compact.
Under this assumption, 
if the steepest descent flow (\ref{eq:a:01}) admits a time-global solution $x(t)$,
then for any sequence $\{t_j\}$ such that $t_j \to \infty$,
the set $\{x(t_j)\}$ is pre-compact. 
Therefore, there exists a subsequence $\{t_{j_k}\}$ such that $\{x(t_{j_k})\}$ converges to some $x_\infty \in X$.
Moreover, if the limit $x_\infty$ is unique, then $x(t)$ converges to $x_\infty$.
\begin{example}[harmonic oscillator]
  \label{example:a:1}
  Consider a periodic harmonic oscillator consisting of three points.
  Let $x_1$, $x_2$, $x_3 \in \R$ be the coordinates of these points,
  and consider a system where they are periodically connected by springs.
  Since this system is invariant under translations 
  $(x_1, x_2, x_3) \mapsto (x_1 + \alpha, x_2 + \alpha, x_3 + \alpha)$ 
  for any $\alpha \in \R$, 
  we may fix $x_2 = 0$ and assume $x_1 < 0 < x_3$. 
  Thus, the system is characterized by the variables $(x_1, x_3) \in \R^2$.
  The energy $E \colon \R^2 \to \R_{\ge0}$ is defined by 
  \begin{equation*}
    \label{eq:a:1:1}
    E(x_1,  x_3)
    = 
    \frac{1}{2}
    \left(
      |x_1|^2
      +
      |x_3|^2
      +
      |x_3 - x_1 - b|^2
    \right), 
  \end{equation*}
  where $b > 0$ represents the period.
  Since $E$ is bounded from below, we consider the steepest descent flow 
  \begin{equation}
    \label{eq:a:1:2}
    \frac{d\v{x}}{dt} = 
    -\nabla E(\v{x}), 
    \quad \text{where } \v{x} = (x_1, x_3)^T, 
  \end{equation}
  and a solution $\v{x}(t)$ satisfies the energy inequality
  \begin{equation*}
    \label{eq:a:1:3}
    E(\v{x}(t)) \le E(\v{x}(0)).
  \end{equation*}
  By computing $\nabla E$, we can rewrite (\ref{eq:a:1:2}) as 
  \begin{equation*}
    \label{eq:a:1:4}
    \frac{d}{dt}
    \begin{bmatrix}
      x_1 \\ x_3
    \end{bmatrix}
    =
    \begin{bmatrix}
      -2 & 1 \\
      1 & -2 \\
    \end{bmatrix}
    \begin{bmatrix}
      x_1 \\ x_3
    \end{bmatrix}
    + 
    \begin{bmatrix}
      -b \\ b
    \end{bmatrix}, 
  \end{equation*}
  i.e.,
  \begin{equation}
    \label{eq:a:1:5}
    \frac{d\v{x}}{dt} = A \v{x} + \v{b}.
  \end{equation}
  It is easy to obtain the solution of (\ref{eq:a:1:5}) 
  \begin{equation}
    \label{eq:a:1:6}
    \v{x}(t) = e^{tA} \v{x}(0) + A^{-1} e^{tA} \v{b} - A^{-1} \v{b}, 
  \end{equation}
  and the solution extends to $t \to \infty$.
  Since $A$ is negative definite, the solution (\ref{eq:a:1:6}) satisfies
  \begin{equation*}
    \label{eq:a:1:7}
    \lim_{t \to \infty} \v{x}(t) 
    =
    - A^{-1} \v{b}
  \end{equation*}
  for any initial value $\v{x}(0)$, 
  implying the existence of a unique stable solution.
\end{example}
\begin{example}[Particles with Coulomb repulsive force]
  \label{example:a:2}
  Consider two particles with like charges (repulsive force).
  Let $z_1, z_2 \in \R$ be the coordinates, with an interaction potential $|z_1 - z_2|^{-1}$.
  Fixing $z_2 = 0$ and assuming $z_1 > 0$, the energy is $E(z_1) = 1/z_1$. 
  While $E$ is bounded from below, 
  the sub-level sets $L(c) = \{ x \in \R \mid x > 0, E(x) \le c\}$ are non-compact (as $x \to \infty$).
  The steepest descent flow is 
  \begin{equation*}
    \label{eq:a:2:1}
    E(z_1)
    = 
    \frac{1}{z_1}, 
    \quad
    z_1 > 0.
  \end{equation*}
  Since the energy functional $E$ is bounded from below, 
  we may consider the steepest descent flow 
  \begin{equation}
    \label{eq:a:2:2}
    \frac{dz_1}{dt} = 
    -\nabla E(z_1), 
  \end{equation}
  and a solution $z_1(t)$ satisfies the energy inequality
  \begin{equation*}
    \label{eq:a:2:3}
    E(z_1(t)) \le E(z_1(0)).
  \end{equation*}
  However, any sub-level set $L(c) = \{ x \in \R | x > 0, E(x) \le c\}$ is non-compact in $\R$.
  By computing $\nabla E$, we may rewrite (\ref{eq:a:2:2}) as 
  \begin{equation}
    \label{eq:a:2:4}
    \frac{dz_1}{dt}
    =
    \frac{1}{z_1^2}
  \end{equation}
  It is easy to obtain the solution of (\ref{eq:a:2:4}) 
  \begin{equation}
    \label{eq:a:2:5}
    z_1(t) = (3 t + z_1(0)^3)^{1/3}
  \end{equation}
  and the solution extends to $t \to \infty$, 
  but the solution (\ref{eq:a:2:5}) diverges to infinity for any initial value with $O(t^{1/3})$.
\end{example}

\begin{remark}
  The energy in the main text is a combination of these two types of energies.
  The $O(t^{1/3})$ growth rate in Theorem \ref{claim:weaving:ode:untangle} arises from the same mechanism as in Example \ref{example:a:2}.
\end{remark}

\begin{remark}
  There are several approach to identify stable configurations of entangled structure by introducing energy of knots
  (cf.~\cite{ohara, simon}).
  In those study, the structures are considered as $1$-dimensional objects.
  It would be however more useful to consider graphs as $0$-dimensional objects (atoms) with interactions
  between atoms and find stable configurations in terms of energy of $0$-dimensional structures when we apply chemistry
  and materials science.
\end{remark}

\section{Appendix: Pair of identical graphs: the case of square lattices}
\label{sec:appendix:b}
A {\em pair of identical graphs} is a periodic graph on the $xy$-plane and its copy arranged in $\R^3$.
Similar results to those for weaves can be shown for a pair of identical graphs.
Here we illustrate ideas of proof by using a pair of identical square lattices which are entangled alternatively as a basic example.
We describe the essential components required to establish these results for square lattices.
\par
Let $X = (V, E)$ be a square lattice in the $xy$-plane; 
namely, $X$ is periodic with respect to two linearly independent directions.
Analogous to the case of weaves, we consider a pair $(X, S)$, 
where $S \colon V \to \{1, -1\}$ is the {\em crossing information}.
\begin{definition}
  A pair $(X, S)$ is called an {\em pair of identical graphs}, 
  or shortly {\em entangled}, 
  if there exist $u$, $v \in V$ such that $S(u) S(v) = -1$, 
  otherwise we call it {\em untangled}.
\end{definition}
\par
The two graphs are exactly the same, except for the $z$-coordinates of each vertex.
Thus, it should be noted that an entangled pair will always intersect at an edge.
\par
Several approaches have been proposed to identify stable configurations 
of entangled systems through the introduction of knot energies 
(cf.~\cite{ohara, simon}).
In these studies, the structures are treated as one-dimensional objects, 
and repulsive interactions are imposed between knotted curves.
In contrast, we model graphs as zero-dimensional objects, 
vertices, and introduce interactions between them indicated by edges, 
motivated by atomic configurations in chemistry and materials science. 
The energy considered in our framework arises from attractive forces between
adjacent vertices within a graph and repulsive forces 
between corresponding vertices of a pair of graphs.
\par
Let $(X, S)$ be a pair of identical square graphs.
A pair of non-parallel periodic translations $\tau_1$ and $\tau_2$ of $W$ induce linearly independent vectors 
$\iota_1$, $\iota_2 \in \mathbb{Z}\times\mathbb{Z}$, respectively.
Namely, for any vertex $v_{i,j}$ of $X$, we have $\tau_1(c_{i,j})=c_{(i,j) + \iota_1}$ and $\tau_2(c_{i,j})=c_{(i,j) + \iota_2}$.
Set $V = \{ (i, j) \in \Z \times \Z \mid (i, j) = \alpha \iota_1 + \beta \iota_2, 0 \leq \alpha < 1, 0 \leq \beta < 1 \}$.
\par
By the definition of a pair of identical square graphs, 
a square graph consists of vertices $(x_{i,j},y_{i,j},z^r_{i,j})$, 
and another square graph consists of vertices $(x_{i,j},y_{i,j},z^b_{i,j})$.
These determine the following maps:
\begin{align*}
  \vfunc{x}&\colon \Z\times\Z\to\R^2, & 
  z_r&\colon \Z\times\Z\to\R,&
  z_b&\colon \Z\times\Z\to\R;\\
  \vfunc{x}&(i,j)=(x_{i,j},y_{i,j}), & z_r&(i,j)=z^r_{i,j}, & z_b&(i,j)=z^b_{i,j}.
\end{align*}
We call the triple $(\vfunc{x},z_r,z_b)$ a \emph{configuration} of the pair of identical square graphs $(X, S)$.
\par
For any index $u = (i, j) \in \Z \times \Z$, there exist unique real numbers $\alpha$, $\beta$ 
such that $u = \alpha \iota_1 + \beta \iota_2$.
We define the projection $p\colon \Z \times \Z \to V$ and the shift $\vfunc{b}\colon \Z \times \Z \to \R^2$ by
\begin{align*}
  p(i, j) &= (i, j) - \lfloor \alpha \rfloor \iota_1 - \lfloor \beta \rfloor \iota_2,\\
  \vfunc{b}(i, j) &= \lfloor \alpha \rfloor \tau_1 + \lfloor \beta \rfloor \tau_2.
\end{align*}
\par
We define the energy $E(X, S)$ of the pair of identical square graphs $(X, S)$ by
\begin{equation}
  \label{eq:entangled:energy}
  \begin{aligned}
    E(X, S)
    &=
    \frac{1}{2} \sum_{v \in V} \sum_{u \in E_v} \lVert\vfunc{x}(p(u)) + \vfunc{b}(u) - \vfunc{x}(v)\rVert^2
    \\
    &+
    \frac{1}{2} \sum_{v \in V} \sum_{u \in E_v} \lvert z_r(p(u)) - z_r(v) \rvert^2
    \\
    &+
    \frac{1}{2} \sum_{v \in V} \sum_{u \in E_v} \lvert z_b(p(u)) - z_b(v) \rvert^2
    \\
    &+
    \sum_{v \in V} \frac{1}{\lvert z_b(v) - z_r(v) \rvert}.
  \end{aligned}
\end{equation}
Here, we remark that $E_v = \{(i+1, j), \, (i-1, j), \, (i, j+1)\, (i, j-1)\}$, since $X$ is a square graph and $v = (i, j)$.
Throughout this subsection, we may also write $E(X, S)$ as $E(\Phi)$ using the configuration $\Phi=(\vfunc{x},z_r,z_b)$ of $(X, S)$.
\par
\begin{figure}[H]
  \centering
  \begin{tabular}{llll}
    \raisebox{70pt}{(a)}
    &\includegraphics[bb=0 0 59 59,width=80pt]{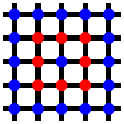}
    &\includegraphics[bb=0 0 278 107,width=80pt]{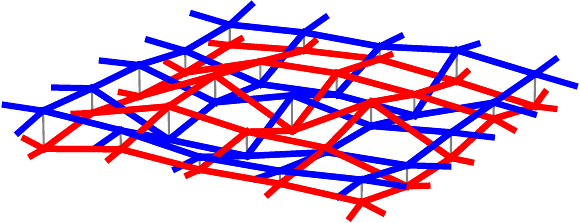}
    &\includegraphics[bb=0 0 278 117,width=80pt]{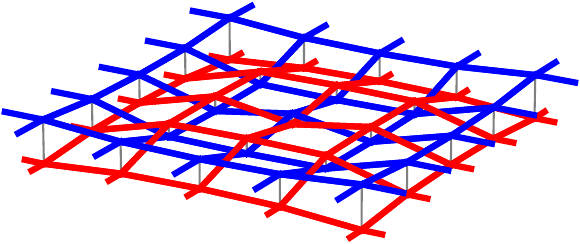}
    \\
    \raisebox{70pt}{(b)}
    &\includegraphics[bb=0 0 59 59,width=80pt]{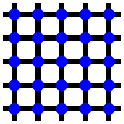}
    &\includegraphics[bb=0 0 278 107,width=80pt]{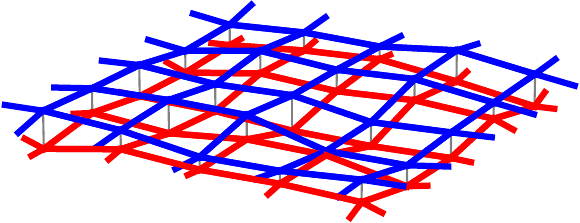}
    &\includegraphics[bb=0 0 280 207,width=80pt]{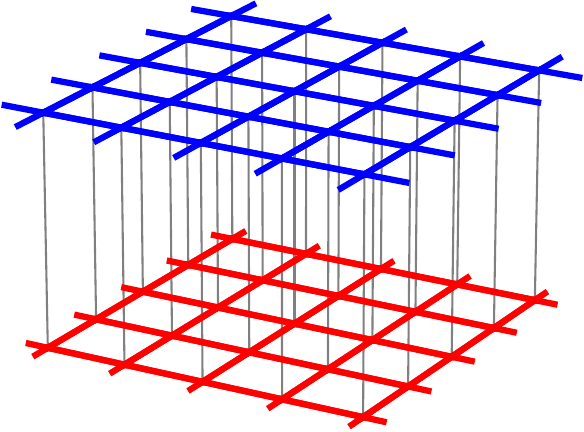}
    \\
  \end{tabular}
  \caption{
    Center column: initial configuration of a pair of identical square graphs in $\R^3$.
    Left column: their projection onto the $xy$-plane.
    Right column: stable configuration according to the energy \eqref{eq:entangled:energy}.
    (a) is entangled, while (b) is untangled.
  }
  \label{fig:entangled:0}
\end{figure}

%
The steepest descent flow of the energy is given as the following system (\ref{eq:entangled:energy}):
\begin{equation}
  \label{eq:entangled:ode:00}
  \begin{aligned}
    \frac{d}{dt} \vfunc{x}(t, v)
    &=
    \sum_{u \in E_v} (\vfunc{x}(t, p(u)) - \vfunc{x}(t, v)) + \sum_{u \in E_v} \vfunc{b}(u), 
    \\
    \frac{d}{dt} z_r(t, v)
    &=
    \sum_{u \in E_v^r}
    (z_r(t, p(u)) - z_r(t, v))
    - \frac{S(v)}{|z_b(t, v) - z_r(t, v)|^2}, 
    \\
    \frac{d}{dt} z_b(t, v)
    &=
    \sum_{u \in E_v^b}
    (z_b(t, p(u)) - z_b(t, v))
    + \frac{S(v)}{|z_b(t, v) - z_r(t, v)|^2}, 
    \\
  \end{aligned}
\end{equation}
The energy \eqref{eq:entangled:energy} and the steepest descent flow \eqref{eq:entangled:ode:00} are analogous to those for weaves.
In particular, the linear parts (the squared difference terms between adjacent points in the energy) are identical to the weave case, except for the specific adjacency relations.
Moreover, the non-linear terms (the terms with negative powers) are also the same as in weaves.
Therefore, we consider the existence of time-global solutions to \eqref{eq:entangled:ode:00} and the stable configurations of the pairs of graphs (see Figure \ref{fig:entangled:0}).
\par
The following proposition provides the key estimate for the global existence of solutions to (\ref{eq:entangled:ode:00}).
\begin{proposition}
  \label{claim:entangled:ode:d}
  If $(X, S)$ is entangled (i.e., there exist $u$, $v \in V$ such that $S(u) S(v) = -1$), 
  then there exists a positive constant $C = C(E_0) > 0$, independent of time $t$, such that 
  \begin{displaymath}
    |\v{d}(t)| \le C
  \end{displaymath}
  along a solution of (\ref{eq:entangled:ode:00}), 
  where $d(t, v) = z_b(t, v) - z_r(t, v)$ and 
  \begin{math}
    \v{d}(t) = 
    \begin{bmatrix}
      d(t, v)
    \end{bmatrix}_{v \in X}.
  \end{math}
\end{proposition}
The proof of Proposition \ref{claim:entangled:ode:d} is nearly identical to that of Proposition \ref{claim:weaving:ode:minimal:d}. 
Consequently, we obtain the following results (cf. Theorems \ref{claim:weaving:ode:1-component} and \ref{claim:weaving:ode:untangle}).
\begin{theorem}
  Assume that a pair of identical square graphs $(X, S)$ is entangled. 
  Then, the maximal interval of existence for the solution to \eqref{eq:entangled:ode:00} is $[0, \infty)$.
  Furthermore, for a suitable sequence $\{t_j\}$ with $t_j \to \infty$,
  the solution $(\v{x}(t), \v{z}_b(t), \v{z}_r(t))$ converges to a unique stationary solution $(\v{X}, \v{Z}_b, \v{Z}_r)$.
\end{theorem}
\begin{theorem}
  If a pair of identical square graphs $(X, S)$ is untangled,
  the two graphs move apart at a rate of $t^{1/3}$ along the steepest descent flow 
  \eqref{eq:entangled:ode:00} of the energy \eqref{eq:entangled:energy}.
\end{theorem}
Here, the symbols $\v{x}$, $\v{z}_r$, and $\v{z}_b$ are vectors consisting of their corresponding functions.

\section{Appendix: Proof of Lemma~\ref{lemma:weaving:height order:partial order}}
\label{sec:appendix:partial-order}

In this appendix, we will show Lemma~\ref{lemma:weaving:height order:partial order}.

\begin{proof}[Proof of Lemma~\ref{lemma:weaving:height order:partial order}]

  We can immediately verify reflexivity by Definition~\ref{definition:weaving:height_order}~(1).
The antisymmetry and transitivity will be proved in the following two claims.

\begin{claim}
  The relation $\preceq$ is antisymmetric.
\end{claim}

\begin{claimproof}

  Let $W_1$ and $W_2$ be weave-connected components of a weave $W$ such that $W_1 \preceq W_2$ and $W_2 \preceq W_1$.
  By Remark~\ref{remark:weaving:height order:either one} and Lemma~\ref{lemma:weaving:height order:comparability}, either $W_1 = W_2$ or both $W_1$ and $W_2$ are untangled threads without any crossings with each other.

  We will show that the latter case cannot occur.
  We sufficiently consider the case where $W_1$ and $W_2$ are both red threads $r_1$ and $r_2$.
  Since $W_1 \preceq W_2$ by the asssumption, there exists a blue thread $b$ such that the crossing between $W_1$ ($W_2$, respectively) and $b$ has a $+1$ sign ($-1$, respectively).
  Similarly, by $W_2 \preceq W_1$, there exists a blue thread $b'$ such that the crossing between $W_2$ ($W_1$, respectively) and $b'$ has a $+1$ sign ($-1$, respectively).
  So, the set $\{r_1, r_2, b, b' \}$ forms an \AltQuad.
  This contradicts the assumption that both $r_1$ and $r_2$ are untangled threads.
  Hence, if $W_1 \preceq W_2$ and $W_2 \preceq W_1$, then $W_1 = W_2$.
\end{claimproof}

\begin{claim}
  The relation $\preceq$ is transitive.
\end{claim}

\begin{claimproof}
  Let $W_1$, $W_2$, and $W_3$ be weave-connected components of a weave $W$ such that $W_1 \preceq W_2$ and $W_2 \preceq W_3$.
  We sufficiently consider the case where $W_1$, $W_2$ and $W_3$ are all distinct.

  We first consider the case where $W_1$ and $W_3$ are parallel threads.
  Here, we assume that they are both red threads $r_1$ and $r_3$.

  If $W_2$ is not an untangle red thread, by $W_1 \preceq W_2$ and $W_2 \preceq W_3$, there exists a blue thread in $W_2$ such that the crossing between $W_1$ ($W_3$, respectively) and the blue thread has a $+1$ sign ($-1$, respectively).
  Hence, we obtain $W_1 \preceq W_3$.

  We assume that $W_2$ is an untangled red thread $r_{2}$.
  By $W_1 \preceq W_2$, there exists a blue thread $b$ such that the crossing between $W_1$ ($W_2$, respectively) and $b$ has $+1$ sign ($-1$, respectively).
  Similarly, by $W_2 \preceq W_3$, there exists a blue thread $b'$ such that the crossing between $W_2$ ($W_3$, respectively) and $b'$ has a $+1$ sign ($-1$, respectively).
  So, if the crossing sign between $b$ and $r_3$ is $+1$, the subweave $\{r_{2}, r_{3}, b, b' \}$ forms an \AltQuad, which contradicts the assumption that $r_{2}$ and $r_{3}$ are untangled threads (see Fig.~\ref{fig:weaving:height_order_transitivity}(a)).
  Hence, the crossing between $b$ and $r_{3}$ has a $-1$ sign, and we obtain $W_1 \preceq W_3$.

  \begin{figure}
    \begin{minipage}[b]{0.48\hsize}
      \centering
      \includegraphics[page=1, width=0.55\linewidth]{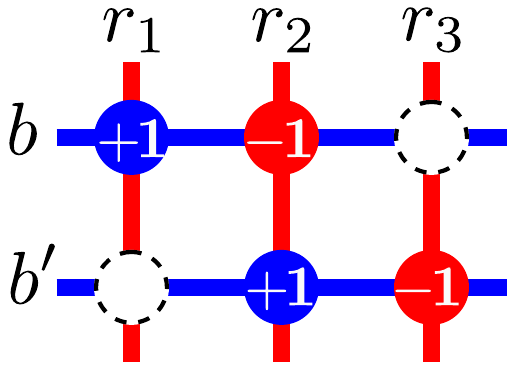}\\
      (a)
    \end{minipage}
    \begin{minipage}[b]{0.48\hsize}
      \centering
      \includegraphics[page=2, width=0.55\linewidth]{Figure/height_order_transitivity.pdf}\\
      (b)
    \end{minipage}
    \caption{Illustration for the proof of Claim~2 in Lemma~\ref{lemma:weaving:height order:partial order}. (a) The case where $W_1$ and $W_3$ are parallel untangled threads. (b) The case where $W_1$ and $W_3$ intersect. This figure only shows the case where $W_2$ contains an \AltQuad.}
    \label{fig:weaving:height_order_transitivity}
  \end{figure}

  We next consider the case where $W_1$ and $W_3$ have an intersection.
  By Lemma~\ref{lemma:weaving:height order:comparability}, either $W_1 \preceq W_3$ or $W_3 \preceq W_1$ holds.
  We assume that $W_3 \preceq W_1$.
  So, if $W_{1}$ and $W_{3}$ contain a red thread $r_{1}$ and a blue thread $b_{3}$, respectively, then the crossing between them has $-1$ sign.
  Since $W_1 \preceq W_2$ and $W_2 \preceq W_3$, there exist a blue thread $b$ and a red thread $r$ such that the following conditions hold (see Fig.~\ref{fig:weaving:height_order_transitivity}(b)):
  \begin{enumerate}
    \item the crossing between $r_{1}$ and $b$ has a $+1$ sign,
    \item the crossing between $b_{3}$ and $r$ has a $+1$ sign, and
    \item the crossings between $b$ and $r$ has a $-1$ sign.
  \end{enumerate}
  Thus, the set $\{r_{1}, b_{3}, b, r\}$ forms an \AltQuad, which contradicts the assumption that both $W_1$ and $W_3$ are distinct weave-connected components.
  If the case where $W_{1}$ and $W_{3}$ contain a blue thread $b_1$ and a red thread $r_3$, respectively, can be proved in the same way.
  Therefore, we obtain $W_1 \preceq W_3$.
\end{claimproof}

Therefore, the relation $\preceq$ is a partial order on the set of weave-connected components and a total order on the set of weave layers.
\end{proof}

\subsection*{Acknowledgements}
M.K. and H.N. acknowledge the JSPS Grant-in-Aid for Scientific Research (B) under Grant No. JP23K25769.
H.N. acknowledges the JSPS Grant-in-Aid for Scientific Research (C) under Grant No. JP24K06710. 
E.S. acknowledges the JSPS Grant-in-Aid for Scientific Research (C) under Grant No. JP24K06859.

\par\vspace{\baselineskip}
\begin{flushleft}
  \begin{tabular}{ll}
    Motoko Kotani:  & {\ttfamily motoko.kotani.d3@tohoku.ac.jp}
    \\
    Hisashi Naito:  & {\ttfamily naito@math.nagoya-u.ac.jp}
    \\
    Naoki Sakata: & {\ttfamily sakata@casis.sakura.ne.jp}
    \\
    Eriko Shinkawa: & {\ttfamily eriko.shinkawa.e8@tohoku.ac.jp}
  \end{tabular}
\end{flushleft}

\begin{thebibliography}{99}
\bibitem{bourgoin}
  M.~O.~Bourgoin,
  {Twisted link theory},
  Algebr.~Geom.~Topol. {\bfseries 8} (2008), 1249--1279.
\bibitem{evans-hyde0}
  T.~Castle, M.~E.~Evans, and S.~T.~Hyde,
  {Entanglement of embedded graphs},
  Prog.~Theor.~Phys.~Suppl., {\bfseries 191} (2011), 235--244.
\bibitem{dechant}
  A.~Dechant, T.~Ohto, Y.~Ito, M.~V.~Makarova, Y.~Kawabe, T.~Agari, H.~Kumai, Y.~Takahashi, H.~Naito, and M.~Kotani,
  {Geometric model of 3D curved graphene with chemical dopants},
  Carbon {\bfseries 182} (2021), 223--232.
\bibitem{diamantis1}
  I.~Diamantis, S.~Lambropoulou, and S.~Mahmoudi,
  {Equivalence of dounbly periodic tangles},
  Mathematics {\bfseries 14 }(2026), 1071.
\bibitem{diamantis2}
  I.~Diamantis, S.~Lambropoulou, and S.~Mahmoudi,
  {Directional Invariants of Doubly Periodic Tangles},
  Symmetry, {\bfseries 16} (2024), 968.
\bibitem{evans-hyde1}	
  M.~E.~Evans, V.~Robins, and S.~T.~Hyde,
  {Periodic entanglement I: networks from hyperbolic reticulations},
  Acta Cryst., {\bfseries A69} (2013), 241--261.
\bibitem{evans-hyde2}
  M.~E.~Evans, V.~Robins, and S.~T.~Hyde,
  {Periodic entanglement II: weavings from hyperbolic line patterns},
  Acta Cryst., {\bfseries A69} (2013), 262--275.
\bibitem{evans-hyde3}
  M.~E.~Evans and S.~T.~Hyde, 
  {Periodic entanglement III: tangled degree-3 finite and layer net intergrowths from rare forests},
  Acta Cryst., {\bfseries A71} (2015), 599--611.
\bibitem{fks1}
  M.~Fukuda, M.~Kotani, and S.~Mahmoudi, 
  {Classification of doubly periodic untwisted $(p,q)$-weaves by their crossing number}, 
  J. Knot Theory Ramif., {\bfseries 32} (2023), 2350032.
\bibitem{fks2}
  M.~Fukuda, M.~Kotani, and S.~Mahmoudi, 
  {Construction of weaving and polycatenanes motifs from periodic tilings of the plane}, 
  J. Knot Theory Ramif., {\bfseries 35} (2026), 2550070.
\bibitem{grishanov1}
  S.~Grishanov, V.~Meshkov, and A.~Omelchenko,
  {A topological study of textile structures. Part I: an introduction to topological methods},
  Text. Res. J., {\bfseries 79}
  (2009), 702--713.
\bibitem{grishanov2}
  S.~Grishanov, V.~Meshkov, and A.~Omelchenko,
  {A topological study of textile structures. Part II: topological invariants in application to textile structures}.
  Text. Res. J., {\bfseries 79} (2009), 822--836.
\bibitem{grishanov3}
  S.~A.~Grishanov and V.~A.~Vassiliev,
  {Invariants of links in 3-manifolds and splitting problem of textile structures},
  J. Knot Theory Ramif., 
  {\bfseries 20}
  (2011), 345--370.
\bibitem{Jost}
  J.~Jost, Riemannian Geometry and Geometric Analysis, 
  Springer, 2017
\bibitem{lkty}
  Y.~Liu,  M.~O'Keeffe, M.~M.~J.~Treacy, and O.~M.~Yaghi, 
  The geometry of periodic knots, polycatenanes and weaving from a chemical perspective: a library for reticular chemistry, 
  Chem.~Soc.~Rev., {\bfseries 47} (2018), 4642--4664.
\bibitem{okeeffe}
  S.~T.~Hyde, B.~Chen, and M.~O'Keeffe,
  {Some equivalent two-dimensional weavings at the molecular scale in 2D and 3D metal organic frameworks}.
  CrystEngComm, 
  {\bfseries 18}
  (2016), 7607--7613.
\bibitem{kawauchi}
  A.~Kawauchi, 
  Complexities of a knitting pattern, 
  React.~Funct.~Polym., 
  {\bfseries 131}
  (2018), 230--236.
\bibitem{Kotani-Sunada}
  M.~Kotani and T.~Sunada, 
  Standard realizations of crystal lattices via harmonic maps, 
  Trans. Amer. Math. Soc., 
  {\bfseries 353}
  (2001) 1--20.
\bibitem{Kotorii-Yoshida}
  Y.~Kotorii and K.~Yoshida, 
  Linking numbers for periodic tangles, 
  arXiv:2509.24568 (2025).
\bibitem{ohara}
  J.~O'Hara, {Energy of a knot}, Topology {\bfseries 30} (1999), 241--247.
\bibitem{Pontryagin}
  L.~S.~Pontryagin, Ordinary Differential Equations, Peragamon Press, 1962.
\bibitem{simon}
  J.~K.~Simon, {Energy functions for polygonal knots}, J.~Knot Theor.~Ramif., 
  {\bfseries 3} (1994), 299--320.
\bibitem{Sunada}
  T.~Sunada, 
  {\em Topological crystallography, With a view towards discrete geometric analysis},
  Springer, 2013.
\end{thebibliography}
\end{document}